\documentclass[11pt,twoside]{article}
\usepackage{amssymb,amsthm}
\usepackage[tbtags]{amsmath}  %
\usepackage{cite}
\usepackage{indentfirst}
\usepackage{cases}
\usepackage{graphicx}
\usepackage{float}
\usepackage{placeins}
\graphicspath{{figure/}}
\usepackage{lineno}
\usepackage{enumerate}
\usepackage{threeparttable}
\usepackage{subcaption}
\usepackage{multirow}
\usepackage{booktabs}
\usepackage{upgreek}  
\usepackage{bm}
\usepackage{stmaryrd}
\usepackage{fancyhdr}
\usepackage{ifthen}
\usepackage{lipsum}
\usepackage{pgfplots}
\pgfplotsset{compat=1.18}

\allowdisplaybreaks[4]
\usepackage[misc]{ifsym}
\usepackage{caption}
\usepackage{appendix}
\usepackage[format=hang]{caption}
\usepackage{algorithm,algorithmic}
\usepackage{charter}  
\usepackage{indentfirst}  
\usepackage{soul}
\usepackage{authblk}

\newcommand{\Rmnum}[1]{\uppercase\expandafter{\romannumeral #1}} 

\usepackage[colorlinks,
linkcolor=blue,
anchorcolor=blue,
citecolor=red]{hyperref}

\numberwithin{equation}{section}

\newtheorem{Lemma}{Lemma}[section]
\newtheorem{Theorem}{Theorem}[section]
\newtheorem{Remark}{Remark}[section]

\newcounter{saveeqn}

\makeatletter
\def\@maketitle{%
	\newpage
	\null
	\vskip 2em%
	\begin{center}%
		\let \footnote \thanks
		{\LARGE \@title \par}%
		\vskip 1.5em%
		{\large
			\lineskip .5em%
			\begin{tabular}[t]{c}%
				\@author
			\end{tabular}\par}%
	\end{center}%
	\par
	\vskip 1.5em}
\makeatother

\makeatletter
\evensidemargin
\oddsidemargin
\makeatother
\title{\bf  Pressure-Robust Enriched Galerkin Methods for Navier-Stokes Equations}
\author{
	Chun Song\footnote{School of Mathematics, Sichuan University, Chengdu, Sichuan 610064, China (songchun20001023@gmail.com).}, 
	\ and Minfu Feng\footnote{School of Mathematics, Sichuan University, Chengdu, Sichuan 610064, China (fmf@scu.edu.cn). The work of this author was supported by the National Natural Science Foundation of China (Grant No. 11971337).}
}

\begin{document}
	\maketitle
	\newcommand\blfootnote[1]{%
		\begingroup
		\renewcommand\thefootnote{}\footnote{#1}%
		\par\setlength\parindent{2em}
		\endgroup
	}
	\captionsetup[figure]{labelfont={bf},labelformat={default},labelsep=period,name={Fig.}}
	\captionsetup[table]{labelfont={bf},labelformat={default},labelsep=period,name={Tab.}}
\begin{abstract}

This paper develops a pressure-robust enriched Galerkin (PR-EG) finite element method for the stationary incompressible Navier--Stokes equations.
Starting from the low-order $P_1$--$P_0$ finite element pair, the velocity space is enriched by an elementwise piecewise linear function, yielding an inf--sup stable velocity--pressure pair with only a minimal increase in the number of degrees of freedom.
A symmetric interior penalty formulation is employed for the viscous term, while a reconstruction operator $\mathcal{R}$ is incorporated into the right-hand side and the nonlinear convection term to achieve pressure robustness.
As a result, the velocity error estimate is independent of the pressure approximation, thereby reducing the influence of gradient forces on the discrete velocity and improving the robustness of the method in pressure-dominated and small-viscosity regimes.
Under suitable regularity and small-data assumptions, the well-posedness of the discrete problem is established and a priori error estimates are derived.
Numerical experiments confirm the predicted convergence behavior and, in particular, demonstrate the theoretically expected pressure-robust performance of the proposed method in high-Reynolds-number flows and in problems dominated by gradient forces.

\


\end{abstract}

	\baselineskip 15pt
	\parskip 10pt
	\setcounter{page}{1}
	\vspace{-0.5cm}
\section{Introduction}

The stationary incompressible Navier--Stokes equations constitute one of the fundamental mathematical models for viscous fluid flow and arise in a wide range of applications in fluid dynamics, geophysical modeling, energy systems, and biomedical flows \cite{batchelor2000introduction,anderson2011ebook,vallis2017atmospheric,masset2020linear,colladoslara2022datadriven,rotavera2021influence,wedi2020baseline,overholt2022carbon,roby2021neurovascular}.
Their numerical approximation remains challenging because of the simultaneous presence of the incompressibility constraint and the nonlinear convective term.
These difficulties become more pronounced as the viscosity decreases and the relative importance of convection increases.
At the analytical level, the nonlinear structure of the Navier--Stokes equations also leads to fundamental questions concerning existence, uniqueness, and regularity \cite{fefferman2006existence}.
Consequently, the construction of stable and accurate discretizations for incompressible flow continues to be an important topic in numerical analysis.\\
\indent For mixed finite element discretizations, the approximation spaces for velocity and pressure must satisfy an appropriate inf--sup condition in order to guarantee stability of the saddle-point problem \cite{boffi2013mixed}.
Classical conforming finite element methods provide efficient discretizations but do not, in general, simultaneously provide local conservation and strong control of convection-dominated flow.
Discontinuous Galerkin (DG) methods offer greater flexibility in the treatment of interelement fluxes and are particularly attractive for locally conservative and convection-dominated problems \cite{di2011mathematical,john2017divergence,linke2008divergence,tezduyar1983finite,muralikrishnan2020multilevel,dang2019numerical}.
This flexibility is accompanied by a larger number of degrees of freedom and a more complicated algebraic structure, and balancing stability, accuracy, and computational complexity remains an important consideration in the design of numerical schemes \cite{binkowski2020riesz,jiang2019linearly}.
Enriched Galerkin (EG) methods provide an intermediate framework between continuous and discontinuous Galerkin discretizations.
By enriching a continuous finite element space with a small number of discontinuous functions, EG methods retain many of the local features of DG methods while using substantially fewer degrees of freedom \cite{yi2022enriched,lee2024uniform}.\\
\indent For incompressible flow problems, Yi et al.\ developed an EG discretization for the Stokes equations using continuous piecewise linear velocities enriched by one discontinuous piecewise linear component per element together with a piecewise constant pressure space \cite{yi2022enriched}.
Lee and Mu subsequently introduced a pressure-robust EG formulation for the Brinkman equations by employing a velocity reconstruction operator into an $H(\mathrm{div},\Omega)$-conforming space \cite{lee2024uniform}.
Related low-order and bubble-stabilized discontinuous Galerkin constructions were investigated by Burman and Stamm for elliptic and Stokes problems \cite{burman2009low,burman2010bubble}.
For the stationary incompressible Navier--Stokes equations, Rivi\`{e}re and Sardar developed and analyzed a low-order DG formulation in which special treatment of the nonlinear upwind term plays an essential role in the stability and convergence analysis \cite{riviere2014penalty}.
Extending the low-order EG framework to the nonlinear Navier--Stokes equations therefore requires not only a stable velocity--pressure pair but also a convection discretization that preserves the appropriate stability properties.\\
\indent A further difficulty arises from the interaction between the pressure and the discrete velocity.
For many inf--sup stable mixed methods, the velocity approximation may be polluted by the discrete approximation of gradient forces.
In particular, pressure-dependent contributions in velocity error estimates can become significant when the viscosity is small.
Pressure-robust methods are designed to eliminate this effect by ensuring that irrotational forces are represented through the pressure and do not generate spurious discrete velocities.
This principle is closely related to the Helmholtz--Hodge decomposition and to the accurate representation of the divergence-free component of the momentum balance \cite{linke2014role,lederer2018hybrid}.
Related approaches for controlling pressure and constraint effects have also been considered in projection-type and structure-preserving discretizations for incompressible and multiphysics problems \cite{shen2020rotational,terekhov2023pressure,eremin2022energy}.\\
\indent In the present work, we develop a pressure-robust enriched Galerkin method for the stationary incompressible Navier--Stokes equations.
We first extend the standard EG formulation to the nonlinear problem by combining a symmetric interior penalty discretization of the viscous term with an upwind treatment of convection.
We then introduce a reconstruction operator $\mathcal{R}$ to enhance the pressure robustness of the velocity approximation.
The reconstructed test function is used in the load term, while the nonlinear form is modified by reconstructing the transported and test variables.
This construction leaves the viscous and velocity--pressure coupling forms unchanged, so that the basic structure of the original EG discretization is preserved.
At the same time, the gradient component of the forcing term is absorbed by the discrete pressure and therefore does not produce a spurious contribution to the discrete velocity.
As a consequence, the velocity approximation is less affected by pressure errors, which is particularly beneficial in pressure-dominated and small-viscosity regimes.
The reconstructed nonlinear form also retains the stability properties required for the analysis of the Navier--Stokes discretization.
Under suitable regularity, we analyze the existence and uniqueness of the discrete solution and derive a priori error estimates.
In particular, the velocity error estimate for the PR-EG formulation contains no pressure-approximation term, in contrast with the corresponding estimate for the standard EG method.\\
\indent The theoretical results are examined through several numerical experiments.
Smooth manufactured solutions are used to study convergence and the behavior of the methods as the viscosity decreases.
A no-flow test with a pure gradient force is used to examine the pressure robustness of the proposed method.
Finally, the lid-driven cavity problem is considered to evaluate the performance of the method in convection-dominated regimes, including computations at high Reynolds numbers.\\
\indent The remainder of the paper is organized as follows.
Section~2 introduces the model problem, notation, mesh-dependent norms, and approximation operators.
Section~3 presents the standard enriched Galerkin discretization of the stationary Navier--Stokes equations.
Section~4 discusses the well-posedness and error analysis of the standard EG method.
In Section~5, the velocity reconstruction operator and the pressure-robust EG formulation are introduced, followed by the corresponding stability and error analysis.
Section~6 presents the numerical experiments for both the standard and pressure-robust schemes.
Finally, Section~7 summarizes the main conclusions.

\section{Model Problem and Preliminaries}

We consider the stationary Navier-Stokes equations on a convex polygon domain $\Omega \subset  \mathbb{R}^d  $ with d = 2,3. Given $\mathbf{f}$ in $\textbf{L}^2(\Omega)$ find $(\mathbf{u},p)$ in $\textbf{H}^1_0(\Omega) \times L^2_0(\Omega)$, such that 
\begin{align}
	-\mu \Delta \mathbf{u}+\mathbf{u} \cdot \nabla \mathbf{u}+\nabla p &  = \mathbf{f} \quad \text{in}~ \Omega, \\
	\nabla \cdot \mathbf{u} & = 0 \quad  \text{in}~ \Omega, \\
	\mathbf{u} & = 0 \quad  \text{on}~ \partial \Omega .
\end{align}

Now, we introduce some useful notions, which are crucial in the following sections. 
Let $\mathcal{T}_h$ be a shape-regular triangulation of $\Omega$ with mesh size $h$. 
And $T$ is the element of it, which means that $T \in \mathcal{T}_h$ .
We denote by $\mathcal{E}_h$ the set of all edges of the triangulation $\mathcal{T}_h$, besides $\mathcal{E}_h = \mathcal{E}_h^o \cup \mathcal{E}_h^b$, where $\mathcal{E}_h^o$ denotes all the interior edges and $\mathcal{E}_h^b$ denotes all the boundary edges. 
For any edge $e \in \mathcal{E}_h$, we denote by $T^+$ and $T^-$ the two elements sharing the edge $e$. 
We denote by $\mathbf{n}_e$ the unit normal vector to the edge $e$ pointing from $T^+$ to $T^-$.
For $e \in \mathcal{E}_h^b$, $\mathbf{n}_e$ is the outward normal vector on the boundary $\partial \Omega$.\\
\indent Additionally, the Sobolev space $H^s(\mathcal{D})$ of a bounded Lipschitz domain $\mathcal{D} \in \mathbb{R}^d (d = 2, 3)$ is equipped with the norm and semi-norm $\left \| \cdot \right \| _{s,\mathcal{D}}$ and $\left | \cdot \right | _{s,\mathcal{D}}$, for any $s \in \mathcal{R}, s \ge 0$.
$H^0(\mathcal{D})$is the same as $L^2(\mathcal{D})$. 
Moreover, $(\cdot, \cdot)_{\mathcal{D}}$ is the inner product in $L^2(\mathcal{D})$.
If we have $\mathcal{D} = \Omega$, we omit the subscript $\mathcal{D}$.
As for vector-valued and tensor-valued Sobolev spaces, we have the same notation.
What's more, we denote $H_0^1(\mathcal{D})$ as follows:
\begin{align*}
	H_0^1(\mathcal{D}) = \{ v \in H^1(\mathcal{D}) : v|_{\partial \mathcal{D}} = 0 \}.
\end{align*}
And $L_0^2(\mathcal{D})$ is denoted as:
\begin{align*}
	L_0^2(\mathcal{D}) = \{ v \in L^2(\mathcal{D}) :(v, 1)_{\mathcal{D}} = 0 \}.
\end{align*}
$P_k(\mathcal{D})$ is the polynomial spaces of degree less than or equal to $k$ on $\mathcal{D}$.
Finally, we introduce $H(\operatorname{div}, \mathcal{D}):=\left\{\mathbf{v} \in\left[L^{2}(\mathcal{D})\right]^{d}: \operatorname{div} \mathbf{v} \in L^{2}(\mathcal{D})\right\}$, which is equipped with the norm:
\begin{align*}
	\left\| \mathbf{v} \right\|^2_{H(\text{div}, \mathcal{D})} := \left\| \mathbf{v} \right\|^2_{0, \mathcal{D}} + \left\| \text{div} \mathbf{v} \right\|^2_{0, \mathcal{D}}.
\end{align*}
In this part, we define the broken Sobolev space as follows:
\begin{align*}
	H^s(\mathcal{T}_h) := \left\{ v \in L^2(\Omega) : v|_T \in H^s(T), \forall T \in \mathcal{T}_h \right\}.
\end{align*}
And its norm is defined as:
\begin{align*}
	\|v\|_{s,\mathcal{T}_h} := \left( \sum_{T \in \mathcal{T}_h} \|v\|_{s,T}^2 \right)^{1/2}.
\end{align*}
When $s = 0$, the $L^2$ -inner product over the mesh $\mathcal{T}_h$ is represented by $(\cdot,\cdot)_{\mathcal{T}_h}$.
Similarly, the $L^2$ -inner product on the set $\mathcal{E}_h$ is expressed as $\left \langle \cdot, \cdot \right \rangle_{\mathcal{E}_h}$, with the corresponding $L^2$ -norm on $\mathcal{E}_h$ is defined in usual manner:
\begin{align*}
	\|v\|_{0, \mathcal{E}_h} := \left( \sum_{e \in \mathcal{E}_h} \|v\|_{0, e}^2 \right)^{1/2}.
\end{align*}
On the broken Sobolev space, we define the piecewise polynomial space as follows:
\begin{align*}
	P_k(\mathcal{T}_h) = \{ v \in L^2(\Omega) \colon v|_T \in P_k(T), \, \forall T \in \mathcal{T}_h \}.
\end{align*}
Finally, we introduce the jump and average of $v$ on $e \in \mathcal{E}_h$:
\begin{align*}
	[v] := \begin{cases} 
		v^+ - v^- & \text{on } e \in \mathcal{E}_h^o, \\
		\hphantom{v^+ - } v & \text{on } e \in \mathcal{E}_h^b,
	\end{cases} \quad 
	\{ v \} := \begin{cases} 
		(v^+ + v^-)/2 & \text{on } e \in \mathcal{E}_h^o, \\
		\hphantom{v^+ + } v & \text{on } e \in \mathcal{E}_h^b.
	\end{cases}
\end{align*}
$v^+$ and $v^-$ denote the trace of $v|_+$ and $v|_-$ on $e \in \partial{T^+} \cap \partial{T^-}$ respectively.
This definition is also valid in vector-valued functions together with tensor-valued functions.
Now, we let's recall a significant property in $H^1(T)$,
\begin{align}
	\|v\|_{0,e}^2 \leq C \left( h_T^{-1} \|v\|_{0,T}^2 + h_T \|\nabla v\|_{0,T}^2 \right), \quad \forall v \in H^1(T).
\end{align}
We introduce the discrete $H^1$-norm over $[H^1_0(\Omega)]^d$,
\begin{align}
	\| \mathbf{v} \|_{\mathcal{E}}^2 := \| \nabla \mathbf{v} \|_{0, \mathcal{T}_h}^2 + \rho \| h_{e}^{-1/2} [\mathbf{v}] \|_{0, \mathcal{E}_h}^2, \quad \forall \mathbf{v} \in [H^1_0(\Omega)]^d, \label{eq:d_H1_norm}
\end{align}
$\rho$ is positive constant, a penalty parameter. 

\section{Enriched Galerkin Methods for Navier-Stokes Equations}

Following the Enriched Galerkin framework of Lee et al.\cite{lee2024uniform}, we employ the same finite-dimensional velocity space and focus on its extension to the nonlinear stationary Navier--Stokes equations.
\begin{align*}
	&\mathbf{V}_h := \mathbf{C}_h \oplus \mathbf{D}_h,\\
	&\mathbf{C}_h=\{\mathbf{v}_h^C \in[H_0^1(\Omega)]^d : \mathbf{v}_h^C|_T \in[P_1(T)]^d, \forall T \in \mathcal{T}_h\},\\
	&\mathbf{D}_h = \{ \mathbf{v}_h^D \in L^2(\Omega) \colon \mathbf{v}_h^D|_T = c(\mathbf{x} - \mathbf{x}_T), c \in \mathbb{R}, \forall T \in \mathcal{T}_h \}.\\
\end{align*}

$\mathbf{C}_h$ is the continuous space for velocity, and $\mathbf{D}_h$ is the discontinuous space for velocity. Besides, $\mathbf{x}_T$ is the barycenter of $ T \in \mathcal{T}_h$. Thus, we can see that velocity consists of a continuous and discontinuous part, which means that $\mathbf{v}_h = \mathbf{v}_h^C + \mathbf{v}_h^D$. \\
\indent On the other hand, we construct presure space as follows: 
\begin{align*}
	\mathbf{Q}_h := \{ q_h \in L^2_0(\Omega) \colon q_h|_T \in P_0(T),\, \forall T \in \mathcal{T}_h \}.
\end{align*}

Thus, we can obtain the EG methods for Navier-Stokes euqations, which means that we want to find $(\mathbf{u}_h,p_h) \in \mathbf{V}_h \times Q_h$ satisfying the following equations:
\begin{align}
	\mu \mathbf{a}(\mathbf{u}_h, \mathbf{v}_h) - \mathbf{b}(\mathbf{v}_h, p_h) + \mathbf{b}(\mathbf{u}_h, q_h) + \mathbf{c}(\mathbf{u}_h, \mathbf{u}_h, \mathbf{u}_h, \mathbf{v}_h) = (\mathbf{f},\mathbf{v}_h), \quad \forall (\mathbf{v}_h,q_h) \in \mathbf{V}_h \times Q_h . \label{eq:model_problem}
\end{align}
And the above bilinear forms are defined as follows:
\begin{align}
	\mathbf{a}(\mathbf{u}, \mathbf{v}) :=&\, (\nabla \mathbf{u}, \nabla \mathbf{v})_{\mathcal{T}_h} - \langle\{\nabla \mathbf{u}\}\mathbf{n}_e, [\mathbf{v}]\rangle_{\mathcal{E}_h} \nonumber \\
	& - \langle\{\nabla \mathbf{v}\}\mathbf{n}_e, [\mathbf{u}]\rangle_{\mathcal{E}_h} + \rho\langle h_e^{-1}[\mathbf{u}], [\mathbf{v}]\rangle_{\mathcal{E}_h}, \label{eq_variation_a}\\
	\mathbf{b}(\mathbf{v}, q) :=&\, (\nabla \cdot \mathbf{v}, q)_{\mathcal{T}_h} - \langle[\mathbf{v}] \cdot \mathbf{n}_e, \{q\}\rangle_{\mathcal{E}_h}. \label{eq_variation_b}
\end{align}
Besides, the nonlinear convection term is discretized with the notion $\mathbf{v}^{\text{int}}$ and $\mathbf{v}^{\text{ext}}$, which refers to the restriction of $\mathbf{v}$ to the element $T$.
The unit normal outwards vector of $T$ is $\mathbf{n}_T$.
\begin{align}
	\mathbf{c}(\mathbf{z}, \mathbf{u}; \mathbf{v}, \mathbf{w}) =\,& (\mathbf{u} \cdot \nabla \mathbf{v}, \mathbf{w})_{\mathcal{T}_h} + \frac{1}{2}((\nabla \cdot \mathbf{u})\mathbf{v}, \mathbf{w})_{\mathcal{T}_h} - \frac{1}{2}([\mathbf{u}], \{\mathbf{v} \cdot \mathbf{w}\})_{\mathcal{E}_h} \nonumber \\
	& + \sum_{T \in \mathcal{T}_h} \int_{\partial T^{\mathbf{z}}_-} |\{\mathbf{u}\} \cdot \mathbf{n}_T| (\mathbf{v}^{\text{int}} - \mathbf{v}^{\text{ext}}) \cdot \mathbf{w}^{\text{int}}, \quad \forall \mathbf{z}, \mathbf{u}, \mathbf{v}, \mathbf{w} \in \mathbf{V}_h,
\end{align}   
where
\begin{align}
	\partial T^{\mathbf{z}}_- = \{ x \in \partial{T} : \{ \mathbf{z}\} \cdot \mathbf{n}_T < 0\}.
\end{align}

\begin{Lemma}
	There is a constant $C_0$ independent of $h$ such that:
	\begin{align}
		|\mathbf{c}(\mathbf{z}_h, \mathbf{u}_h; \mathbf{v}_h, \mathbf{w}_h)| &\leq C_0 \|\mathbf{u}_h\|_{\mathcal{E}} \|\mathbf{v}_h\|_{\mathcal{E}} \|\mathbf{w}_h\|_{\mathcal{E}}, &&\forall \mathbf{z}_h, \mathbf{u}_h, \mathbf{v}_h, \mathbf{w}_h \in \mathbf{V}_h, \label{eq:c_bound} \\
		\mathbf{c}(\mathbf{u}_h, \mathbf{u}_h; \mathbf{v}_h, \mathbf{v}_h) &\geq 0, &&\forall \mathbf{u}_h, \mathbf{v}_h \in \mathbf{V}_h.
	\end{align}
\end{Lemma}

\section{Well-Posedness of the EG Method and Error Analysis}
We can obtain the coercivity and continuity of $\mathbf{a}(\cdot, \cdot)$ by using the discrete $H^1$-norm \eqref{eq:d_H1_norm}.
\begin{Lemma}
	\begin{align}
		& \mathbf{a}(\mathbf{v}_h, \mathbf{v}_h) \ge \kappa_1 \| \mathbf{v}_h \|_{\mathcal{E}}^2, &&\forall \mathbf{v}_h \in \mathbf{V}_h,\label{eq:lem4.1_1} \\  
		& |\mathbf{a}(\mathbf{u}_h, \mathbf{v}_h)| \le \kappa_2 \| \mathbf{u}_h \|_{\mathcal{E}}\| \mathbf{v}_h \|_{\mathcal{E}}, &&\forall \mathbf{u}_h,\mathbf{v}_h \in \mathbf{V}_h, \label{eq:lem4.1_2}
	\end{align}
where $\kappa_1$ and $\kappa_2$ are positive constants independent of $h$.
\end{Lemma}

We can seperate the errors into numerical and approximation errors, which are denoted by the interpolant operators $\Pi_h$ and $\mathcal{P}_0$.
So that we have:
\begin{align*}
	\boldsymbol{\chi}_h := \mathbf{u} - \Pi_h \mathbf{u}, \quad
	\mathbf{e}_h := \Pi_h \mathbf{u} - \mathbf{u}_h, \quad
	\xi_h := p - \mathcal{P}_0 p, \quad
	\epsilon_h := \mathcal{P}_0 p - p_h.
\end{align*}
We define $\Pi_h: [H^2(\Omega)]^d \rightarrow \mathbf{V}_h$ as follows:
\begin{align}
	\mathbf{\Pi}_h \mathbf{v} = \mathbf{\Pi}_h^C \mathbf{v} + \mathbf{\Pi}_h^D \mathbf{v},
\end{align}
where $\mathbf{\Pi}_h^C \mathbf{v} \in \mathbf{C}_h$ and $\mathbf{\Pi}_h^D \mathbf{v} \in \mathbf{D}_h$, and we have 
\begin{align}
	( \nabla \cdot \Pi_h^D \mathbf{v} , 1 )_T = ( \nabla \cdot ( \mathbf{v} - \Pi_h^C \mathbf{v} ), 1 )_T.
\end{align}
Therefore, we have following bounds:
\begin{align}
	&| \mathbf{v} - \Pi_h \mathbf{v} |_{j, \mathcal{T}_h} \leq C h^{m - j} | \mathbf{v} |_m, \quad &&0 \leq j \leq m \leq 2, \quad \forall \mathbf{v} \in [H^2(\Omega)]^d, \\
	&\| \mathbf{v} - \Pi_h \mathbf{v} \|_{\mathcal{E}} \leq C h \| \mathbf{v} \|_2, \quad &&\forall \mathbf{v} \in [H^2(\Omega)]^d, \label{eq:4_9} \\
	&\| \Pi_h \mathbf{v} \|_{\mathcal{E}} \leq C | \mathbf{v} |_1, \quad &&\forall \mathbf{v} \in [H_0^1(\Omega)]^d. \label{eq:Pi_h}
\end{align}
As for $\mathcal{P}_0$, it is introduced as the $L^2$-projection, $\mathcal{P}_0: H^1(\Omega) \rightarrow Q_h$, which satisfies $(q - \mathcal{P}_0q, 1)_T = 0, \forall T \in \mathcal{T}_h$. It also has the following bound:
\begin{align}
	\| q - \mathcal{P}_0 q \|_0 \leq C h \| q \|_1, \quad \forall q \in H^1(\Omega). \label{eq:q_interpolant}
\end{align}

\subsection{Well-Posedness of the EG Method}
First, let's consider the inf-sup condition of the bilinear form $\mathbf{b}(\cdot, \cdot)$.
\begin{Lemma}
	\label{lem:inf-sup}
	We assume that the penalty parameter $\rho$ is sufficiently large.
	Thus, we can obtain the following inf-sup condition, with the positive constant $C_{IS}$.
	\begin{align}
		\inf_{q_h \in Q_h} \sup_{\mathbf{v}_h \in \mathbf{V}_h} \frac{\mathbf{b}(\mathbf{v}_h, q_h)}{\| \mathbf{v}_h \|_{\mathcal{E}} \| q_h \|_0} \geq C_{IS}.
	\end{align}
\end{Lemma}
\begin{proof}
	From the definition of $\mathbf{V}_h$ and $Q_h$, we can know that $\nabla \cdot \mathbf{V}_h = Q_h$, so there exists $\mathbf{v}_h \in \mathbf{V}_h$ such that $\nabla \cdot \mathbf{v}_h = q_h$ and $| \mathbf{v}_h |_1 \le C^{**}\| q_h \|_0$.
	Besides, we have $( \nabla \cdot ( \mathbf{v}_h - \Pi_h \mathbf{v} ), 1 )_{\mathcal{T}_h} = 0$.
	Thus, by using \eqref{eq:Pi_h} we can derive the following result:
	\begin{align*}
		\frac{b(\Pi_h \mathbf{v}, q_h)}{\|\Pi_h \mathbf{v}\|_{\mathcal{E}}} &= \frac{(\nabla \cdot \Pi_h \mathbf{v}, q_h)_{\mathcal{T}_h}}{\|\Pi_h \mathbf{v}\|_{\mathcal{E}}} - \frac{\left \langle  [\Pi_h \mathbf{v}] \cdot \mathbf{n}_e, \{q_h\}\right \rangle_{\mathcal{E}_h}}{\|\Pi_h \mathbf{v}\|_{\mathcal{E}}} \\
		&=\frac{(\nabla \cdot \mathbf{v}_h, q_h)_{\mathcal{T}_h}}{\|\Pi_h \mathbf{v}\|_{\mathcal{E}}} - \frac{\left \langle  [\Pi_h \mathbf{v}] , \{q_h\}\cdot \mathbf{n}_e\right \rangle_{\mathcal{E}_h}}{\|\Pi_h \mathbf{v}\|_{\mathcal{E}}} \\
		& \ge \frac{\| q_h\|_0^2}{C|\mathbf{v}_h|_1} - \frac{\| \{ q_h \} \cdot \mathbf{n}_e \|_{0,\mathcal{E}_h} \| [\Pi_h \mathbf{v}] \|_{0,\mathcal{E}_h}}{\|\Pi_h \mathbf{v}\|_{\mathcal{E}}}\\
		& \ge \frac{\| q_h\|_0^2}{CC^{**} \| q_h\|_0} - \frac{C(h^{-1/2}\| q_h \|_{0} ) \cdot h^{1/2}( \| \nabla \mathbf{v}_h \|_{0,\mathcal{T}_h}^2 + \rho \| h^{-1/2} [\Pi_h \mathbf{v}] \|_{0,\mathcal{E}_h}^2)^{1/2}}{\|\Pi_h \mathbf{v}\|_{\mathcal{E}}}\\
		& \ge \frac{\| q_h\|_0}{CC^{**}} -C \frac{\| q_h \|_{0} \| \Pi_h \mathbf{v} \|_{\mathcal{E}}}{\|\Pi_h \mathbf{v}\|_{\mathcal{E}}}\\
		& \ge C_{IS} \| q_h\|_0.
	\end{align*}
\end{proof}
Next, we consider the existence of the solution. 
To obtain this result, we define the map $G$ by linearizing the convection term, and prove the existence of the fixed point of $G$.
\begin{align*}
	G: \mathbf{V}_h \times Q_h &\rightarrow \mathbf{V}_h \times Q_h, \\
	(\tilde{\mathbf{u}}_h, \tilde{p}_h) &\mapsto (\mathbf{u}_h, p_h), \\
\end{align*}
where $(\mathbf{u}_h, p_h)$ satisfies the following equation:
\begin{align}
	\label{eq:linearized equation}
	\mu \mathbf{a}(\mathbf{u}_h, \mathbf{v}_h) - \mathbf{b}(\mathbf{v}_h, p_h) + \mathbf{b}(\mathbf{u}_h, q_h) + \mathbf{c}(\tilde{\mathbf{u}}_h, \tilde{\mathbf{u}}_h; \mathbf{u}_h, \mathbf{v}_h) = (\mathbf{f}, \mathbf{v}_h), \quad \forall (\mathbf{v}_h, q_h) \in \mathbf{V}_h \times Q_h.
\end{align}

We denote the left side as $ \mathcal{S}_{\tilde{\mathbf{u}}_h}((\mathbf{u}_h, p_h),(\mathbf{v}_h, q_h))$, and the right side remains the same.
Therefore, we can write \eqref{eq:linearized equation} as $ \mathcal{S}_{\tilde{\mathbf{u}}_h}((\mathbf{u}_h, p_h),(\mathbf{v}_h, q_h)) = (\mathbf{f}, \mathbf{v}_h)$.

\begin{Remark}
	We define the map $G$ for EG method, and it remains the same for the PR-EG method.
\end{Remark}

Assume that $ \| \tilde{\mathbf{u}}_h \|_{\mathcal{E}} \le C $, and define
\begin{align*}
	\tilde{\mathbf{a}}(\mathbf{u}_h, \mathbf{v}_h)
	:=
	\mu \mathbf{a}(\mathbf{u}_h, \mathbf{v}_h)
	+
	\mathbf{c}(\tilde{\mathbf{u}}_h, \tilde{\mathbf{u}}_h; \mathbf{u}_h, \mathbf{v}_h).
\end{align*}
It follows from \eqref{eq:lem4.1_2} and \eqref{eq:c_bound} that
\begin{align*}
	| \tilde{\mathbf{a}}(\mathbf{u}_h, \mathbf{v}_h) |
	&\le
	\kappa_2 \| \mathbf{u}_h \|_{\mathcal{E}} \| \mathbf{v}_h \|_{\mathcal{E}}
	+
	C_N \| \tilde{\mathbf{u}}_h \|_{\mathcal{E}}
	\| \mathbf{u}_h \|_{\mathcal{E}}
	\| \mathbf{v}_h \|_{\mathcal{E}}
	\\
	&\le
	\kappa_2 \| \mathbf{u}_h \|_{\mathcal{E}} \| \mathbf{v}_h \|_{\mathcal{E}}
	+
	C_N C \| \mathbf{u}_h \|_{\mathcal{E}} \| \mathbf{v}_h \|_{\mathcal{E}}
	\\
	&\le
	\tilde{C} \| \mathbf{u}_h \|_{\mathcal{E}} \| \mathbf{v}_h \|_{\mathcal{E}},
	\quad
	\forall \mathbf{u}_h, \mathbf{v}_h \in \mathbf{V}_h,
	\\
	\tilde{\mathbf{a}}(\mathbf{v}_h, \mathbf{v}_h)
	&\ge
	\kappa_1 \| \mathbf{v}_h \|_{\mathcal{E}}^2,
	\quad
	\forall \mathbf{v}_h \in \mathbf{V}_h,
\end{align*}
where $ \tilde{C} = \kappa_2 + C_0 C $.\\
\indent Therefore, the bilinear form $ \tilde{\mathbf{a}}(\cdot, \cdot) $ is continuous and coercive.
Combining these properties with the previously established inf-sup condition, the corresponding linear problem satisfies the standard stability framework.
Consequently, the linearized problem is well posed and admits a unique solution.
Thus, for any given $ (\tilde{\mathbf{u}}_h, \tilde{p}_h) $, there exists a unique solution $ (\mathbf{u}_h, p_h) $.\\
\indent Moreover, the Cauchy--Schwarz inequality gives
\begin{align*}
	| ( \mathbf{f}, \mathbf{v}_h ) |
	\leq
	C_{CS}
	\| \mathbf{f} \|_{L^2(\Omega)}
	\| \mathbf{v}_h \|_{\mathcal{E}},
	\quad
	\forall ( \mathbf{v}_h, q_h ) \in \mathbf{V}_h \times Q_h,
\end{align*}
where $ C_{CS} $ is a positive constant depending only on the domain $ \Omega $.
Hence, the right-hand side is continuous with respect to the energy norm.\\
\indent Taking $ \mathbf{v}_h = \mathbf{u}_h $ and $ q_h = p_h $ in \eqref{eq:linearized equation}, we obtain
\begin{align}
	\mu \| \mathbf{u}_h \|_{\mathcal{E}}^2
	\le
	( \mathbf{f}, \mathbf{u}_h )
	\leq
	C_{CS}
	\| \mathbf{f} \|_{L^2(\Omega)}
	\| \mathbf{u}_h \|_{\mathcal{E}}.
	\label{eq:u_h bound}
\end{align}
Therefore,
\begin{align*}
	\| \mathbf{u}_h \|_{\mathcal{E}}
	\leq
	C_1.
\end{align*}
This estimate provides the boundedness required for constructing a closed bounded set associated with the fixed-point mapping.\\
\indent Define
\begin{align*}
	\mathbf{X}_h
	:=
	\{\mathbf{v}_h \in \mathbf{V}_h : \| \mathbf{v}_h \|_{\mathcal{E}} \leq C_1\}.
\end{align*}
The set $ \mathbf{X}_h $ is closed and bounded in $ \mathbf{V}_h $ and is relatively compact since $ \mathbf{V}_h $ is finite-dimensional.
Restricting the mapping $ G $ to $ \mathbf{X}_h \times Q_h $, the preceding estimate shows that its image remains in this set.
Therefore, $ G $ is well defined.
To apply the Leray--Schauder fixed-point theorem, it remains to verify the continuity of $ G $.\\
\indent Let $ (\tilde{\mathbf{u}}_h^n) \subset \mathbf{X}_h $ satisfy
\begin{align*}
	\tilde{\mathbf{u}}_h^n
	\rightarrow
	\tilde{\mathbf{u}}_h
	\quad
	\text{in } \mathbf{X}_h.
\end{align*}
Similarly, let $ (\tilde{p}_h^n) \subset Q_h $ satisfy
\begin{align*}
	\tilde{p}_h^n
	\rightarrow
	\tilde{p}_h
	\quad
	\text{in } Q_h.
\end{align*}
Denote
\begin{align*}
	G(\tilde{\mathbf{u}}_h^n, \tilde{p}_h^n)
	=
	(\mathbf{u}_h^n, p_h^n),
	\qquad
	G(\tilde{\mathbf{u}}_h, \tilde{p}_h)
	=
	(\mathbf{u}_h, p_h).
\end{align*}
By the definition of the linearized problem, the corresponding solutions satisfy
\begin{align}
	&\mu \mathbf{a}(\mathbf{u}_h-\mathbf{u}_h^n, \mathbf{v}_h)
	-
	\mathbf{b}(\mathbf{v}_h, p_h-p_h^n)
	+
	\mathbf{b}(\mathbf{u}_h-\mathbf{u}_h^n, q_h)
	\nonumber
	\\
	&\quad
	+
	\mathbf{c}(\tilde{\mathbf{u}}_h, \tilde{\mathbf{u}}_h; \mathbf{u}_h, \mathbf{v}_h)
	-
	\mathbf{c}(\tilde{\mathbf{u}}_h^n, \tilde{\mathbf{u}}_h^n; \mathbf{u}_h^n, \mathbf{v}_h)
	=
	0.
	\label{eq:G_difference}
\end{align}

The form $ \mathbf{c} $ can be decomposed into a linear part and a nonlinear part as follows:
\begin{align}
	\mathbf{c}(\mathbf{z}_h, \mathbf{u}_h; \mathbf{v}_h, \mathbf{w}_h)
	=
	\mathbf{c}_{\ell}(\mathbf{u}_h; \mathbf{v}_h, \mathbf{w}_h)
	+
	\mathbf{c}_{n \ell}(\mathbf{z}_h, \mathbf{u}_h; \mathbf{v}_h, \mathbf{w}_h),
	\quad
	\forall \mathbf{u}_h, \mathbf{v}_h, \mathbf{w}_h, \mathbf{z}_h \in \mathbf{V}_h.
\end{align}
Here,
\begin{align*}
	&\mathbf{c}_{\ell}(\mathbf{u}_h; \mathbf{v}_h, \mathbf{w}_h)
	=
	(\mathbf{u}_h \cdot \nabla \mathbf{v}_h, \mathbf{w}_h)_{\mathcal{T}_h}
	+
	\frac{1}{2}((\nabla \cdot \mathbf{u}_h)\mathbf{v}_h, \mathbf{w}_h)_{\mathcal{T}_h}
	-
	\frac{1}{2}([\mathbf{u}_h], \{\mathbf{v}_h \cdot \mathbf{w}_h\})_{\mathcal{E}_h},
	\\
	&\mathbf{c}_{n \ell}(\mathbf{z}_h, \mathbf{u}_h; \mathbf{v}_h, \mathbf{w}_h)
	=
	\sum_{T \in \mathcal{T}_h}
	\int_{\partial T^{\mathbf{z}_h}_{-}}
	|\{\mathbf{u}_h\} \cdot \mathbf{n}_T|
	(\mathbf{v}_h^{\text{int}}-\mathbf{v}_h^{\text{ext}})
	\cdot
	\mathbf{w}_h^{\text{int}},
	\\
	&\hspace{10cm}
	\forall \mathbf{u}_h, \mathbf{v}_h, \mathbf{w}_h, \mathbf{z}_h \in \mathbf{V}_h.
\end{align*}

For the linear part, its bilinear structure gives
\begin{align}
	&\mathbf{c}_{\ell}(\tilde{\mathbf{u}}_h; \mathbf{u}_h, \mathbf{v}_h)
	-
	\mathbf{c}_{\ell}(\tilde{\mathbf{u}}_h^n; \mathbf{u}_h^n, \mathbf{v}_h)
	\nonumber
	\\
	&\quad
	=
	\mathbf{c}_{\ell}(\tilde{\mathbf{u}}_h^n; \mathbf{u}_h-\mathbf{u}_h^n, \mathbf{v}_h)
	+
	\mathbf{c}_{\ell}(\tilde{\mathbf{u}}_h-\tilde{\mathbf{u}}_h^n; \mathbf{u}_h, \mathbf{v}_h),
	\quad
	\forall \mathbf{v}_h \in \mathbf{V}_h.
	\label{eq:linear_difference}
\end{align}

For the nonlinear part, direct expansion yields
\begin{align}
	&\mathbf{c}_{n \ell}(\tilde{\mathbf{u}}_h, \tilde{\mathbf{u}}_h; \mathbf{u}_h, \mathbf{v}_h)
	-
	\mathbf{c}_{n \ell}(\tilde{\mathbf{u}}_h^n, \tilde{\mathbf{u}}_h^n; \mathbf{u}_h^n, \mathbf{v}_h)
	\nonumber
	\\
	&\quad
	=
	\mathbf{c}_{n \ell}
	(\tilde{\mathbf{u}}_h^n, \tilde{\mathbf{u}}_h-\tilde{\mathbf{u}}_h^n;
	\mathbf{u}_h, \mathbf{v}_h)
	\nonumber
	\\
	&\qquad
	+
	\mathbf{c}_{n \ell}
	(\tilde{\mathbf{u}}_h^n, \tilde{\mathbf{u}}_h^n;
	\mathbf{u}_h-\mathbf{u}_h^n, \mathbf{v}_h)
	\nonumber
	\\
	&\qquad
	+
	\mathbf{c}_{n \ell}
	(\tilde{\mathbf{u}}_h, \tilde{\mathbf{u}}_h;
	\mathbf{u}_h, \mathbf{v}_h)
	-
	\mathbf{c}_{n \ell}
	(\tilde{\mathbf{u}}_h^n, \tilde{\mathbf{u}}_h;
	\mathbf{u}_h, \mathbf{v}_h),
	\nonumber
	\\
	&\hspace{10cm}
	\forall \mathbf{v}_h \in \mathbf{V}_h.
	\label{eq:nonlinear_difference}
\end{align}

Substituting the above decompositions into the difference equation gives
\begin{align}
	&\mathcal{S}_{\tilde{\mathbf{u}}_h^n}
	((\mathbf{u}_h-\mathbf{u}_h^n, p_h-p_h^n), (\mathbf{v}_h, q_h))
	\nonumber
	\\
	&\quad
	=
	-\mathbf{c}_{\ell}
	(\tilde{\mathbf{u}}_h-\tilde{\mathbf{u}}_h^n; \mathbf{u}_h, \mathbf{v}_h)
	\nonumber
	\\
	&\qquad
	-
	\mathbf{c}_{n \ell}
	(\tilde{\mathbf{u}}_h^n, \tilde{\mathbf{u}}_h-\tilde{\mathbf{u}}_h^n;
	\mathbf{u}_h, \mathbf{v}_h)
	\nonumber
	\\
	&\qquad
	-
	\mathbf{c}_{n \ell}
	(\tilde{\mathbf{u}}_h, \tilde{\mathbf{u}}_h;
	\mathbf{u}_h, \mathbf{v}_h)
	+
	\mathbf{c}_{n \ell}
	(\tilde{\mathbf{u}}_h^n, \tilde{\mathbf{u}}_h;
	\mathbf{u}_h, \mathbf{v}_h).
	\label{eq:S_difference}
\end{align}

By the continuity estimate \eqref{eq:c_bound}, we have
\begin{align}
	|\mathbf{c}_{\ell}
	(\tilde{\mathbf{u}}_h-\tilde{\mathbf{u}}_h^n;
	\mathbf{u}_h, \mathbf{v}_h)|
	&\leq
	C_N
	\| \tilde{\mathbf{u}}_h-\tilde{\mathbf{u}}_h^n \|_{\mathcal{E}}
	\| \mathbf{u}_h \|_{\mathcal{E}}
	\| \mathbf{v}_h \|_{\mathcal{E}},
	\label{eq:lbound}
	\\
	|\mathbf{c}_{n \ell}
	(\tilde{\mathbf{u}}_h^n, \tilde{\mathbf{u}}_h-\tilde{\mathbf{u}}_h^n;
	\mathbf{u}_h, \mathbf{v}_h)|
	&\leq
	C_N
	\| \tilde{\mathbf{u}}_h-\tilde{\mathbf{u}}_h^n \|_{\mathcal{E}}
	\| \mathbf{u}_h \|_{\mathcal{E}}
	\| \mathbf{v}_h \|_{\mathcal{E}}.
	\label{eq:nlbound}
\end{align}

On the other hand, Lemma \ref{lem:inf-sup} gives the equivalent stability estimate
\begin{align}
	C_{IS}
	\| ( \mathbf{u}_h, p_h ) \|_{\mathbf{V}_h \times Q_h}
	\leq
	\sup_{(\mathbf{v}_h, q_h) \in \mathbf{V}_h \times Q_h}
	\frac{
	\mathcal{S}_{\tilde{\mathbf{u}}_h}
	((\mathbf{u}_h, p_h), (\mathbf{v}_h, q_h))
	}{
	\| ( \mathbf{v}_h, q_h ) \|_{\mathbf{V}_h \times Q_h}
	}.
	\label{eq:inf-sup2}
\end{align}
Here,
\begin{align}
	\| ( \mathbf{v}_h, q_h ) \|_{\mathbf{V}_h \times Q_h}
	:=
	\left(
	\| \mathbf{v}_h \|_{\mathcal{E}}^2
	+
	\| q_h \|_0^2
	\right)^{1/2}.
	\label{vq_norm}
\end{align}

Combining \eqref{eq:lbound}, \eqref{eq:nlbound}, and \eqref{eq:inf-sup2}, we obtain
\begin{align}
	&C_{IS}
	\| (\mathbf{u}_h-\mathbf{u}_h^n, p_h-p_h^n) \|_{\mathbf{V}_h \times Q_h}
	\nonumber
	\\
	&\quad
	\leq
	C_N
	\left(
	\| \tilde{\mathbf{u}}_h \|_{\mathcal{E}}
	+
	\| \mathbf{u}_h \|_{\mathcal{E}}
	\right)
	\| \tilde{\mathbf{u}}_h-\tilde{\mathbf{u}}_h^n \|_{\mathcal{E}}
	\nonumber
	\\
	&\qquad
	+
	\sup_{\mathbf{v}_h \in \mathbf{V}_h}
	\frac{
	|
	\mathbf{c}_{n \ell}
	(\tilde{\mathbf{u}}_h^n, \tilde{\mathbf{u}}_h;
	\mathbf{u}_h, \mathbf{v}_h)
	-
	\mathbf{c}_{n \ell}
	(\tilde{\mathbf{u}}_h, \tilde{\mathbf{u}}_h;
	\mathbf{u}_h, \mathbf{v}_h)
	|
	}{
	\| \mathbf{v}_h \|_{\mathcal{E}}
	}.
	\label{eq:G_continuity_estimate}
\end{align}
This estimate controls the difference between the corresponding solutions in terms of the difference between the prescribed velocities and provides the basis for proving the continuity of $ G $.\\
\indent Thus, by using the following lemma(see Lemma 6.1 in \cite{riviere2014penalty}), we can conclude that $G$ is continuous.
\begin{Lemma}
	Assume that $(\tilde{\mathbf{u}}_h^n) \in \mathbf{X}_h$ converges to $\tilde{\mathbf{u}}_h \in \mathbf{X}_h$, we have 
	\begin{align}
		\lim_{n\rightarrow\infty} \sup_{\mathbf{v}_h\in \mathbf{V}_h} \frac{|\mathbf{c}_{n\ell}(\tilde{\mathbf{u}}_h^n, \tilde{\mathbf{u}}_h, \mathbf{u}_h, \mathbf{v}_h) - \mathbf{c}_{n\ell}(\tilde{\mathbf{u}}_h, \tilde{\mathbf{u}}_h, \mathbf{u}_h, \mathbf{v}_h)|}{\| \mathbf{v}_h \|_{\mathcal{E}}} = 0, \qquad \forall \mathbf{v}_h \in \mathbf{V}_h.
	\end{align}
\end{Lemma}
We can easily check the compactness of the map $G$, as it can maps a bounded set to a bounded set.
Above all, we can conclude our result by using the Leray-Schauder theorem.
In order to establish the following theorem, we first present the following lemma (see Proposition 4.10 in \cite{girault2009dg}).
\begin{Lemma}
	\begin{align}
		\label{eq:lem4.5}
		\left| \mathbf{c}_{n\ell}(\mathbf{z}_h, \mathbf{w}_h, \mathbf{z}_h, \mathbf{v}_h) - \mathbf{c}_{n\ell}(\mathbf{w}_h, \mathbf{w}_h, \mathbf{z}_h, \mathbf{v}_h) \right|
		\le C \| \mathbf{z}_h - \mathbf{w}_h \|_{\mathcal{E}} \| \mathbf{v}_h \|_{\mathcal{E}} \| \mathbf{z}_h \|_{\mathcal{E}},
		\quad \forall \mathbf{z}_h, \mathbf{w}_h, \mathbf{v}_h \in \mathbf{V}_h ,
	\end{align}	
	where $C$ is a positive constant.
\end{Lemma}
\begin{Theorem}
	The model problem \eqref{eq:model_problem} has a unique solution $(\mathbf{u}_h, p_h) \in \mathbf{V}_h \times Q_h$, under the small data assumption $\frac{C_N}{\mu^2}\| \mathbf{f} \|_{L^2(\Omega)} < 1$,
	where $C_N$ is a constant.
\end{Theorem}

\begin{proof}
	We can choose two solutions of \eqref{eq:model_problem}, $(\mathbf{u}_h^1, p_h^1)$ and $(\mathbf{u}_h^2, p_h^2)$.
	So, we have 
	\begin{align}
		&\mu \mathbf{a}(\mathbf{u}_h^1 -\mathbf{u}_h^2, \mathbf{v}_h) - \mathbf{b}(\mathbf{v}_h, p_h^1 - p_h^2) + \mathbf{b}(\mathbf{u}_h^1 - \mathbf{u}_h^2, q_h) + \mathbf{c}(\mathbf{u}_h^1, \mathbf{u}_h^1, \mathbf{u}_h^1, \mathbf{v}_h) - \mathbf{c}(\mathbf{u}_h^2, \mathbf{u}_h^2, \mathbf{u}_h^2, \mathbf{v}_h) = 0. \label{eq:0_bound}
	\end{align}
	Besides, the left side of \eqref{eq:0_bound} can be written as 
	\begin{align}
		\label{eq:Uniqueness}
		&\mathbf{c}(\mathbf{u}_h^1, \mathbf{u}_h^1, \mathbf{u}_h^1, \mathbf{v}_h) - \mathbf{c}(\mathbf{u}_h^2, \mathbf{u}_h^2, \mathbf{u}_h^2, \mathbf{v}_h) \nonumber\\
		& = \mathbf{c}(\mathbf{u}_h^2, \mathbf{u}_h^2, \mathbf{u}_h^1 - \mathbf{u}_h^2, \mathbf{v}_h) + \mathbf{c}(\mathbf{u}_h^1, \mathbf{u}_h^1 - \mathbf{u}_h^2, \mathbf{u}_h^1, \mathbf{v}_h) + \mathbf{c}_{n \ell}(\mathbf{u}_h^1, \mathbf{u}_h^2, \mathbf{u}_h^1, \mathbf{v}_h) -\mathbf{c}_{n \ell}(\mathbf{u}_h^2, \mathbf{u}_h^2, \mathbf{u}_h^1, \mathbf{v}_h).
	\end{align}
	For convenience, we set $\mathbf{r}_h = \mathbf{u}_h^1 - \mathbf{u}_h^2$.
	By taking $\mathbf{v}_h = \mathbf{r}_h$ and $q_h = p_h^1 - p_h^2$, we have 
	\begin{align*}
		& \mu \| \mathbf{r}_h \|_{\mathcal{E}}^2 + \mathbf{c}(\mathbf{u}_h^1, \mathbf{u}_h^1, \mathbf{u}_h^1, \mathbf{r}_h) - \mathbf{c}(\mathbf{u}_h^2, \mathbf{u}_h^2, \mathbf{u}_h^2, \mathbf{r}_h) = 0.
	\end{align*}
	Combining with \eqref{eq:u_h bound}, \eqref{eq:lem4.5} and \eqref{eq:Uniqueness}, we have 
	\begin{align*}
		\mu \| \mathbf{r}_h \|_{\mathcal{E}}^2 & = - (\mathbf{c}(\mathbf{u}_h^1, \mathbf{r}_h, \mathbf{u}_h^1, \mathbf{r}_h) + \mathbf{c}_{n \ell}(\mathbf{u}_h^1, \mathbf{u}_h^2, \mathbf{u}_h^1, \mathbf{r}_h) -\mathbf{c}_{n \ell}(\mathbf{u}_h^2, \mathbf{u}_h^2, \mathbf{u}_h^1, \mathbf{r}_h)) \\
		& \le | \mathbf{c}(\mathbf{u}_h^1, \mathbf{r}_h, \mathbf{u}_h^1, \mathbf{r}_h) | + | \mathbf{c}_{n \ell}(\mathbf{u}_h^1, \mathbf{u}_h^2, \mathbf{u}_h^1, \mathbf{r}_h) -\mathbf{c}_{n \ell}(\mathbf{u}_h^2, \mathbf{u}_h^2, \mathbf{u}_h^1, \mathbf{r}_h) | \\
		& \le N \| \mathbf{r}_h \|_{\mathcal{E}}^2 \| \mathbf{u}_h^1 \|_{\mathcal{E}} + C \| \mathbf{r}_h \|_{\mathcal{E}}^2 \| \mathbf{u}_h^1 \|_{\mathcal{E}} \\
		& \le \frac{C_N}{\mu} \| \mathbf{r}_h \|_{\mathcal{E}}^2 \| \mathbf{f} \|_{L^2(\Omega)},
	\end{align*}
	where $C_N = max\{ C_P N, C_P C\}$. 
	Thus, it holds that 
	\begin{align*}
		\mu (1 - \frac{C_N}{\mu^2}\| \mathbf{f} \|_{L^2(\Omega)})\| \mathbf{r}_h \|_{\mathcal{E}}^2 \le 0
	\end{align*}
	From the small data assumption, we can see that $\mathbf{r}_h = 0$.
	Through the inf-sup condition, we can also know that $p_h^1 - p_h^2 = 0$.
\end{proof}

\subsection{Error Analysis of the EG Method}

\begin{Lemma}
	$\forall \mathbf{v}_h \in \mathbf{V}_h, q_h \in Q_h$, we have 
	\begin{align}
		\| \mathbf{e}_h \|_{\mathcal{E}}
		&\leq
		Ch
		\left(
		\mu \| \mathbf{u} \|_2
		+
		\| p \|_1
		\right),
		\\
		\| \epsilon_h \|_0
		&\leq
		Ch
		\left(
		\mu \| \mathbf{u} \|_2
		+
		\| p \|_1
		\right).
	\end{align} 
\end{Lemma}
\begin{proof}
	We first derive the error equation between the continuous solution and the discrete solution. Subtracting the discrete formulation from the continuous formulation gives
	\begin{align}
		&\mu \mathbf{a}(\mathbf{e}_h,\mathbf{v}_h)
		-\mathbf{b}(\mathbf{v}_h,\epsilon_h)
		+\mathbf{b}(\mathbf{e}_h,q_h)
		+\mathbf{c}(\mathbf{u},\mathbf{u};\mathbf{u},\mathbf{v}_h)
		-\mathbf{c}(\mathbf{u}_h,\mathbf{u}_h;\mathbf{u}_h,\mathbf{v}_h)
		\nonumber\\
		&\qquad=
		-\mu \mathbf{a}(\boldsymbol{\chi}_h,\mathbf{v}_h)
		+\mathbf{b}(\mathbf{v}_h,\xi_h)
		-\mathbf{b}(\boldsymbol{\chi}_h,q_h).
	\end{align}

	Since $\mathbf{u}\in [H_0^1(\Omega)]^d$, we obtain
	\begin{align}
		&\mathbf{c}(\mathbf{u},\mathbf{u};\mathbf{u},\mathbf{v}_h)
		-\mathbf{c}(\mathbf{u}_h,\mathbf{u}_h;\mathbf{u}_h,\mathbf{v}_h)
		\nonumber\\
		&\qquad=
		\mathbf{c}(\mathbf{u}_h,\mathbf{u};\mathbf{u},\mathbf{v}_h)
		-\mathbf{c}(\mathbf{u}_h,\mathbf{u}_h;\mathbf{u}_h,\mathbf{v}_h)
		\nonumber\\
		&\qquad=
		\mathbf{c}(\mathbf{u}_h,\mathbf{u}_h;\mathbf{e}_h,\mathbf{v}_h)
		+\mathbf{c}(\mathbf{u}_h,\mathbf{e}_h;\Pi_h\mathbf{u},\mathbf{v}_h)
		\nonumber\\
		&\qquad\quad
		+\mathbf{c}(\mathbf{u}_h,\boldsymbol{\chi}_h;
		\Pi_h\mathbf{u},\mathbf{v}_h)
		+\mathbf{c}(\mathbf{u}_h,\mathbf{u};
		\boldsymbol{\chi}_h,\mathbf{v}_h).
	\end{align}

	On the other hand, the inf-sup condition \eqref{eq:inf-sup2} gives
	\begin{align}
		C_{IS}
		\|(\mathbf{e}_h,\epsilon_h)\|_{\mathbf{V}_h\times Q_h}
		&\leq
		\sup_{(\mathbf{v}_h,q_h)\in\mathbf{V}_h\times Q_h}
		\frac{
		\mathcal{S}_{\mathbf{u}_h}
		((\mathbf{e}_h,\epsilon_h),(\mathbf{v}_h,q_h))
		}{
		\|(\mathbf{v}_h,q_h)\|_{\mathbf{V}_h\times Q_h}
		} \\ \nonumber
		&\le \sup_{(\mathbf{v}_h,q_h)\in\mathbf{V}_h\times Q_h} \frac{1}{\|(\mathbf{v}_h,q_h)\|_{\mathbf{V}_h\times Q_h}}(-\mu \mathbf{a}(\boldsymbol{\chi}_h,\mathbf{v}_h) + \mathbf{b}(\mathbf{v}_h,\xi_h) - \mathbf{b}(\boldsymbol{\chi}_h,q_h)\\ \nonumber
		&-\mathbf{c}(\mathbf{u}_h,\mathbf{e}_h;\Pi_h\mathbf{u},\mathbf{v}_h)
		-\mathbf{c}(\mathbf{u}_h,\boldsymbol{\chi}_h;
		\Pi_h\mathbf{u},\mathbf{v}_h)
		-\mathbf{c}(\mathbf{u}_h,\mathbf{u};
		\boldsymbol{\chi}_h,\mathbf{v}_h))
	\end{align}

	Furthermore, the consistency property of the interpolation operator gives
	\begin{align}
		\mathbf{b}(\mathbf{v}-\Pi_h\mathbf{v},q_h)
		=
		0,
		\quad
		\forall q_h\in Q_h.
	\end{align}

	By the continuity of the nonlinear form and the stability of the interpolation operator, we have
	\begin{align}
		\left|
		\mathbf{c}(\mathbf{u}_h,\mathbf{e}_h;
		\Pi_h\mathbf{u},\mathbf{v}_h)
		\right|
		&\leq
		C
		\|\mathbf{e}_h\|_{\mathcal{E}}
		\|\Pi_h\mathbf{u}\|_{\mathcal{E}}
		\|\mathbf{v}_h\|_{\mathcal{E}}
		\nonumber\\
		&\leq
		C
		\|\mathbf{e}_h\|_{\mathcal{E}}
		\|\mathbf{u}\|_1
		\|\mathbf{v}_h\|_{\mathcal{E}}.
	\end{align}

	Moreover, the remaining interpolation error terms can be bounded uniformly as
	\begin{align}
		&\left|
		-\mu \mathbf{a}(\boldsymbol{\chi}_h,\mathbf{v}_h)
		+\mathbf{b}(\mathbf{v}_h,\xi_h)
		-\mathbf{c}(\mathbf{u}_h,\boldsymbol{\chi}_h;
		\Pi_h\mathbf{u},\mathbf{v}_h)
		-\mathbf{c}(\mathbf{u}_h,\mathbf{u};
		\boldsymbol{\chi}_h,\mathbf{v}_h)
		\right|
		\nonumber\\
		&\qquad\leq
		\left| \mu \mathbf{a}(\boldsymbol{\chi}_h,\mathbf{v}_h) \right|
		+\left| \mathbf{b}(\mathbf{v}_h,\xi_h) \right|
		+\left| \mathbf{c}(\mathbf{u}_h,\boldsymbol{\chi}_h;
		\Pi_h\mathbf{u},\mathbf{v}_h) \right|
		+\left| \mathbf{c}(\mathbf{u}_h,\mathbf{u};
		\boldsymbol{\chi}_h,\mathbf{v}_h) \right|
		\nonumber\\
		&\qquad\leq
		Ch
		\left(
		\mu\|\mathbf{u}\|_2+\|p\|_1
		\right)
		\|\mathbf{v}_h\|_{\mathcal{E}}.
	\end{align}

	Combining the above estimates with the stability inequality yields
	\begin{align}
		\|(\mathbf{e}_h,\epsilon_h)\|_{\mathbf{V}_h\times Q_h}
		\leq
		Ch
		\left(
		\mu\|\mathbf{u}\|_2+\|p\|_1 + 
		\right)
		+ C
		\|\mathbf{e}_h\|_{\mathcal{E}}
		\|\mathbf{u}\|_1.
	\end{align}
	Finally, by the definition of the norm
	$\| (\mathbf{e}_h,\epsilon_h) \|_{\mathbf{V}_h \times Q_h}$,
	the velocity and pressure error estimates can be obtained separately as
	\begin{align}
		\| \mathbf{e}_h \|_{\mathcal{E}}
		&\leq
		Ch
		\left(
		\mu \| \mathbf{u} \|_2
		+
		\| p \|_1
		\right),
		\nonumber\\
		\| \epsilon_h \|_0
		&\leq
		Ch
		\left(
		\mu \| \mathbf{u} \|_2
		+
		\| p \|_1
		\right).
	\end{align}
\end{proof}

From the analysis above, we can derive the error estimates for our problem.
\begin{Theorem}
	We assume that $(\mathbf{u}, p) \in \textbf{H}^1_0(\Omega) \times L^2_0(\Omega)$ is the solution of Navier-Stokes equations and $(\mathbf{u}_h, p_h)$ is the solution of \eqref{eq:model_problem} with small data assumption, we have the following error estimates:
	\begin{align}
		&\| \mathbf{u}-\mathbf{u}_h \|_{\mathcal{E}}
		\leq
		Ch
		\left(
		(\mu+1)\| \mathbf{u} \|_2
		+
		\| p \|_1
		\right), \\
		&\| p-p_h \|_0
		\leq
		Ch
		\left(
		\mu\| \mathbf{u} \|_2
		+
		\| p \|_1
		\right).
	\end{align}
\end{Theorem}
\begin{proof}
	We first decompose the total error into the interpolation error and the discrete error. By the triangle inequality, the velocity error satisfies
	\begin{align*}
		\| \mathbf{u}-\mathbf{u}_h \|_{\mathcal{E}}
		\leq
		\| \boldsymbol{\chi}_h \|_{\mathcal{E}}
		+
		\| \mathbf{e}_h \|_{\mathcal{E}}.
	\end{align*}

	For the interpolation error, the approximation property of $\Pi_h$ gives
	\begin{align*}
		\| \boldsymbol{\chi}_h \|_{\mathcal{E}}
		\leq
		Ch\| \mathbf{u} \|_2.
	\end{align*}
	Moreover, the discrete error estimate yields
	\begin{align*}
		\| \mathbf{e}_h \|_{\mathcal{E}}
		\leq
		Ch
		\left(
		\mu\| \mathbf{u} \|_2
		+
		\| p \|_1
		\right).
	\end{align*}
	Combining the above estimates, we obtain
	\begin{align}
		\| \mathbf{u}-\mathbf{u}_h \|_{\mathcal{E}}
		\leq
		Ch
		\left(
		(\mu+1)\| \mathbf{u} \|_2
		+
		\| p \|_1
		\right).
	\end{align}

	Similarly, we can obtain the error estimate for pressure in the same way
	\begin{align}
		\| p-p_h \|_0
		\leq
		Ch
		\left(
		\mu\| \mathbf{u} \|_2
		+
		\| p \|_1
		\right).
	\end{align}
\end{proof}

\section{Pressure-Robust Enriched Galerkin Methods for Navier-Stokes Equations}
We can derive a pressure-robust EG method for Navier-Stokes equations with the velocity reconstruction operator \cite{lee2024uniform}, which is defined as: $\mathcal{R}: \mathbf{V}_h \to \mathcal{BDM}_1(\mathcal{T}_h) \subset H(div,\Omega)$, such that
\begin{align}
&\ \int_e (\mathcal{R} \mathbf{v}) \cdot \mathbf{n}_e p_1 \, ds = \int_e \{\mathbf{v}\} \cdot \mathbf{n}_e p_1 \, ds, & \forall p_1 \in P_1(e), \, \forall e \in \mathcal{E}_h^o,\\
&\ \int_e (\mathcal{R}\mathbf{v}) \cdot \mathbf{n}_e p_1 \, ds = 0, & \forall p_1 \in P_1(e), \, \forall e \in \mathcal{E}_h^b,
\end{align}
$\mathcal{BDM}_1(\mathcal{T}_h)$ is the Brezzi-Douglas-Marini space of degree 1 over $\mathcal{T}_h$.
Combining the EG method and the velocity reconstuction operator, we can obtain the pressure-robust EG method (PR-EG).
\begin{align}
    \label{eq:variation euqation_PREG}
    \mu \mathbf{a}(\mathbf{u}_h, \mathbf{v}_h) - \mathbf{b}(\mathbf{v}_h, p_h) + \mathbf{b}(\mathbf{u}_h, q_h) + \tilde{\mathbf{c}}(\mathbf{u}_h, \mathbf{u}_h; \mathbf{u}_h, \mathbf{v}_h) = (\mathbf{f},\mathcal{R} \mathbf{v}), \quad \forall (\mathbf{v}_h,q_h) \in \mathbf{V}_h \times Q_h.
\end{align}
$\mathbf{a}(\mathbf{u}_h, \mathbf{v}_h)$ and $\mathbf{b}(\mathbf{v}_h, q_h)$ remain the same with \eqref{eq_variation_a} and \eqref{eq_variation_b} respectively. And $\tilde{\mathbf{c}}(\mathbf{u}_h, \mathbf{u}_h, \mathbf{u}_h, \mathbf{v}_h)$ is defined as follows:
\begin{align}
    \tilde{\mathbf{c}}(\mathbf{u}_h, \mathbf{u}_h; \mathbf{u}_h, \mathbf{v}_h) =& \mathbf{c}(\mathbf{u}_h, \mathbf{u}_h; \mathcal{R}\mathbf{u}_h, \mathcal{R}\mathbf{v}_h)\\
    = &(\mathbf{u}_h \cdot \nabla \mathcal{R}\mathbf{u}_h, \mathcal{R}\mathbf{v}_h)_{\mathcal{T}_h} + \frac{1}{2}((\nabla \cdot \mathbf{u}_h)\mathcal{R}\mathbf{u}_h, \mathcal{R}\mathbf{v}_h)_{\mathcal{T}_h} \nonumber \\
    & -\frac{1}{2}([\mathbf{u}_h], \{\mathcal{R}\mathbf{u}_h \cdot \mathcal{R}\mathbf{v}_h\})_{\mathcal{E}_h} + \sum_{T \in \mathcal{T}_h} \int_{\partial T^{\mathbf{u}_h}_-} |\{\mathbf{u}_h\} \cdot \mathbf{n}_T| (\mathcal{R}\mathbf{u}^{\text{int}}_h - \mathcal{R}\mathbf{u}^{\text{ext}}_h) \cdot \mathcal{R}\mathbf{v}_h^{\text{int}}. 
\end{align}

For the reconstrcution operator $\mathcal{R}$, we have the following lemma
\begin{Lemma}
    \label{lem:lem4.2}
    \begin{align}
        \| \mathcal{R} \mathbf{v}_h - \mathbf{v}_h \|_0 \leq C h \| \mathbf{v}_h \|_{\mathcal{E}}, \quad \forall \mathbf{v}_h \in \mathbf{V}_h.
    \end{align}
\end{Lemma}
\begin{proof}
For any $\phi_k \in RT_0(K)$, by using scaling argument, we have the following bound:
\begin{align}
    \| \boldsymbol{\phi}_K \|_{0,K} \leq C \sum_{e \in \partial K} h_e^{1/2} \| \boldsymbol{\phi}_K \cdot \mathbf{n}_K \|_{0,e} \label{eq:scaling argument}
\end{align}
where $C$ is positive constant independent of $h$.
Thus, for $\boldsymbol{\phi}_K = (\mathcal{R}\mathbf{v}_h^D - \mathbf{v}_h^D)|_K \in RT_0(K)$, and $\{ \mathbf{v}_h^D\} - \mathbf{v}_h^D = \pm \frac{1}{2}[\mathbf{v}_h^D]$, we have 
\begin{align}
    \| \mathcal{R} \mathbf{v}_h - \mathbf{v}_h \|_0^2 &= \| \mathcal{R} \mathbf{v}_h^D - \mathbf{v}_h^D \|_{0, \mathcal{T}_h}^2 \leq C \| h_e^{1/2} [\mathbf{v}_h^D] \cdot \mathbf{n}_e \|_{0, \mathcal{E}_h}^2 \leq C h^2 \| h_e^{-1/2} [\mathbf{v}_h^D] \|_{0, \mathcal{E}_h}^2 \nonumber\\
    &= C h^2 \| h_e^{-1/2} [\mathbf{v}_h] \|_{0, \mathcal{E}_h}^2 \leq C h^2 \| \mathbf{v}_h \|_{\mathcal{E}}^2.
\end{align}
\end{proof}
\begin{Lemma}
Under the same assumptions as in Lemma\ref{lem:lem4.2}, there exists a positive
constant $C$, independent of $h$, such that
\[
\|\mathcal{R}\mathbf{v}_h\|_{1}
\leq
C\|\mathbf{v}_h\|_{\mathcal{E}},
\qquad
\forall\,\mathbf{v}_h\in\mathbf{V}_h.
\]
\end{Lemma}

\begin{proof}
By the triangle inequality, we have
\begin{align}
\|\nabla\mathcal{R}\mathbf{v}_h\|_{0}
&\leq
\|\nabla\mathbf{v}_h\|_{0}
+
\|\nabla(\mathcal{R}\mathbf{v}_h-\mathbf{v}_h)\|_{0}.
\end{align}
Since $\mathcal{R}\mathbf{v}_h-\mathbf{v}_h$ is piecewise linear,
the inverse inequality and Lemma\ref{lem:lem4.2} imply
\begin{align}
\|\nabla(\mathcal{R}\mathbf{v}_h-\mathbf{v}_h)\|_{0}
&\leq
Ch^{-1}
\|\mathcal{R}\mathbf{v}_h-\mathbf{v}_h\|_0
\nonumber\\
&\leq
C\|\mathbf{v}_h\|_{\mathcal{E}}.
\end{align}
Therefore,
\begin{align}
\|\nabla\mathcal{R}\mathbf{v}_h\|_{0}
\leq
C\|\mathbf{v}_h\|_{\mathcal{E}}.
\end{align}

Moreover, by the triangle inequality, Lemma\ref{lem:lem4.2}, and the discrete
Poincar\'e inequality, we obtain
\begin{align}
\|\mathcal{R}\mathbf{v}_h\|_0
&\leq
\|\mathbf{v}_h\|_0
+
\|\mathcal{R}\mathbf{v}_h-\mathbf{v}_h\|_0
\nonumber\\
&\leq
C\|\mathbf{v}_h\|_{\mathcal{E}}
+
Ch\|\mathbf{v}_h\|_{\mathcal{E}}
\nonumber\\
&\leq
C\|\mathbf{v}_h\|_{\mathcal{E}}.
\end{align}
Combining the above estimates yields
\[
\|\mathcal{R}\mathbf{v}_h\|_{1}
\leq
C\|\mathbf{v}_h\|_{\mathcal{E}}.
\]
\end{proof}
\begin{Lemma}
There exists a positive constant $C$, independent of $h$, such that
\[
\|\mathcal{R}\mathbf{v}_h\|_{\mathcal{E}}
\leq
C\|\mathbf{v}_h\|_{\mathcal{E}},
\qquad
\forall\,\mathbf{v}_h\in\mathbf{V}_h.
\]
\end{Lemma}

\begin{proof}
For each $K\in\mathcal{T}_h$, set
\[
\boldsymbol{\phi}_K
:=
\left.
\left(
\mathcal{R}\mathbf{v}_h^D-\mathbf{v}_h^D
\right)
\right|_K
\in RT_0(K)\subset[P_1(K)]^d.
\]
Using the same scaling argument as in Lemma~\ref{lem:lem4.2} and the
definition of the reconstruction operator, we obtain
\[
\|\boldsymbol{\phi}_K\|_{0,K}^2
\leq
C
\sum_{e\subset\partial K}
h_e
\|[\mathbf{v}_h^D]\|_{0,e}^2.
\]
Since the mesh is shape-regular,
\[
\sum_{K\in\mathcal{T}_h}
h_K^{-2}\|\boldsymbol{\phi}_K\|_{0,K}^2
\leq
C
\|h_e^{-1/2}[\mathbf{v}_h^D]\|_{0,\mathcal{E}_h}^2.
\]
Moreover, since $\mathbf{v}_h^C$ is continuous and vanishes on
$\partial\Omega$,
\[
[\mathbf{v}_h]=[\mathbf{v}_h^D]
\qquad\text{on }\mathcal{E}_h.
\]
Thus,
\[
\sum_{K\in\mathcal{T}_h}
h_K^{-2}\|\boldsymbol{\phi}_K\|_{0,K}^2
\leq
C
\|h_e^{-1/2}[\mathbf{v}_h]\|_{0,\mathcal{E}_h}^2.
\]

Since $\boldsymbol{\phi}_K\in[P_1(K)]^d$, the inverse inequality gives
\[
\|\nabla\boldsymbol{\phi}_K\|_{0,K}
\leq
Ch_K^{-1}\|\boldsymbol{\phi}_K\|_{0,K}.
\]
Therefore,
\begin{align*}
\|\nabla\mathcal{R}\mathbf{v}_h\|_{0,\mathcal{T}_h}
&\leq
\|\nabla\mathbf{v}_h\|_{0,\mathcal{T}_h}
+
\|\nabla(\mathcal{R}\mathbf{v}_h^D-\mathbf{v}_h^D)\|_{0,\mathcal{T}_h}
\\
&\leq
\|\nabla\mathbf{v}_h\|_{0,\mathcal{T}_h}
+
C\|h_e^{-1/2}[\mathbf{v}_h]\|_{0,\mathcal{E}_h}
\\
&\leq
C\|\mathbf{v}_h\|_{\varepsilon}.
\end{align*}

Next, by the inverse trace inequality,
\[
\|\boldsymbol{\phi}_K\|_{0,e}^2
\leq
Ch_K^{-1}\|\boldsymbol{\phi}_K\|_{0,K}^2,
\qquad
e\subset\partial K.
\]
Since
\[
[\mathcal{R}\mathbf{v}_h]
=
[\mathbf{v}_h]
+
[\mathcal{R}\mathbf{v}_h^D-\mathbf{v}_h^D],
\]
the above estimate yields
\[
\|h_e^{-1/2}[\mathcal{R}\mathbf{v}_h]\|_{0,\mathcal{E}_h}
\leq
C\|h_e^{-1/2}[\mathbf{v}_h]\|_{0,\mathcal{E}_h}.
\]

Therefore, we obtain
\begin{align*}
\|\mathcal{R}\mathbf{v}_h\|_{\varepsilon}^2
&=
\|\nabla\mathcal{R}\mathbf{v}_h\|_{0,\mathcal{T}_h}^2
+
\rho
\|h_e^{-1/2}[\mathcal{R}\mathbf{v}_h]\|_{0,\mathcal{E}_h}^2
\\
&\leq
C\|\mathbf{v}_h\|_{\varepsilon}^2.
\end{align*}
\end{proof}
\begin{Lemma}
	There is a constant $C_0^{PR}$ independent of $h$ such that:
	\begin{align}
		|\tilde{\mathbf{c}}(\mathbf{z}_h, \mathbf{u}_h; \mathbf{v}_h, \mathbf{w}_h)| &\leq C_0^{PR} \|\mathbf{u}_h\|_{\mathcal{E}} \|\mathbf{v}_h\|_{\mathcal{E}} \|\mathbf{w}_h\|_{\mathcal{E}}, &&\forall \mathbf{z}_h, \mathbf{u}_h, \mathbf{v}_h, \mathbf{w}_h \in \mathbf{V}_h, \label{eq:pr_c_bound}\\
		\tilde{\mathbf{c}}(\mathbf{u}_h, \mathbf{u}_h; \mathbf{v}_h, \mathbf{v}_h) &\geq 0, &&\forall \mathbf{u}_h, \mathbf{v}_h \in \mathbf{V}_h.
	\end{align}
\end{Lemma}
\begin{proof}
We first estimate the volume term. By the Sobolev embedding inequality
and the stability property of the reconstruction operator, we obtain
\begin{align}
&
\left|
(\mathbf{u}_h\cdot\nabla \mathcal{R}\mathbf{v}_h,
\mathcal{R}\mathbf{w}_h)_{\mathcal{T}_h}
\right|
\nonumber\\
&\leq
\|\mathbf{u}_h\|_{0,4}
\|\nabla\mathcal{R}\mathbf{v}_h\|_{0}
\|\mathcal{R}\mathbf{w}_h\|_{0,4}
\nonumber\\
&\leq
C
\|\mathbf{u}_h\|_{\mathcal{E}}
\|\mathcal{R}\mathbf{v}_h\|_{1}
\|\mathcal{R}\mathbf{w}_h\|_{1}
\nonumber\\
&\leq
C
\|\mathbf{u}_h\|_{\mathcal{E}}
\|\mathbf{v}_h\|_{\mathcal{E}}
\|\mathbf{w}_h\|_{\mathcal{E}}.
\label{eq:PR_volume1}
\end{align}

Similarly, for the divergence term, we have
\begin{align}
&
\left|
((\nabla\cdot\mathbf{u}_h)\mathcal{R}\mathbf{v}_h,
\mathcal{R}\mathbf{w}_h)_{\mathcal{T}_h}
\right|
\nonumber\\
&\leq
C
\|\nabla\cdot\mathbf{u}_h\|_{0}
\|\mathcal{R}\mathbf{v}_h\|_{0,4}
\|\mathcal{R}\mathbf{w}_h\|_{0,4}
\nonumber\\
&\leq
C
\|\mathbf{u}_h\|_{\mathcal{E}}
\|\mathbf{v}_h\|_{\mathcal{E}}
\|\mathbf{w}_h\|_{\mathcal{E}}.
\label{eq:PR_volume2}
\end{align}

For the jump term, using the trace inequality and the stability
of $\mathcal{R}$, we have
\begin{align}
&
\left|
\left\langle
[\mathbf{u}_h],
\{\mathcal{R}\mathbf{v}_h\cdot\mathcal{R}\mathbf{w}_h\}
\right\rangle_{\mathcal{E}_h}
\right|
\nonumber\\
&\leq
C
\|h_e^{-1/2}[\mathbf{u}_h]\|_{0,\mathcal{E}_h}
\left\|
h_e^{1/2}
\{\mathcal{R}\mathbf{v}_h\cdot\mathcal{R}\mathbf{w}_h\}
\right\|_{0,\mathcal{E}_h}
\nonumber\\
&\leq
C
\|\mathbf{u}_h\|_{\mathcal{E}}
\|\mathcal{R}\mathbf{v}_h\|_{1}
\|\mathcal{R}\mathbf{w}_h\|_{1}
\nonumber\\
&\leq
C
\|\mathbf{u}_h\|_{\mathcal{E}}
\|\mathbf{v}_h\|_{\mathcal{E}}
\|\mathbf{w}_h\|_{\mathcal{E}}.
\label{eq:PR_jump}
\end{align}

The upwind term can be bounded in the same way. Namely,
by the trace inequality,
\begin{align}
&
\left|
\sum_{T\in\mathcal{T}_h}
\int_{\partial T_-^{\mathbf{u}_h}}
|\{\mathbf{u}_h\}\cdot\mathbf{n}_T|
\left(
(\mathcal{R}\mathbf{v}_h)^{\rm int}
-
(\mathcal{R}\mathbf{v}_h)^{\rm ext}
\right)
\cdot
(\mathcal{R}\mathbf{w}_h)^{\rm int}
\right|
\nonumber\\
&\leq
C
\|\mathbf{u}_h\|_{\mathcal{E}}
\|\mathcal{R}\mathbf{v}_h\|_{1}
\|\mathcal{R}\mathbf{w}_h\|_{1}
\nonumber\\
&\leq
C
\|\mathbf{u}_h\|_{\mathcal{E}}
\|\mathbf{v}_h\|_{\mathcal{E}}
\|\mathbf{w}_h\|_{\mathcal{E}}.
\label{eq:PR_upwind}
\end{align}

Combining \eqref{eq:PR_volume1}--\eqref{eq:PR_upwind}, we obtain
\begin{align}
\left|
\tilde{\mathbf{c}}
(\mathbf{z}_h,\mathbf{u}_h;
\mathbf{v}_h,\mathbf{w}_h)
\right|
\leq
C_0^{PR}
\|\mathbf{u}_h\|_{\mathcal{E}}
\|\mathbf{v}_h\|_{\mathcal{E}}
\|\mathbf{w}_h\|_{\mathcal{E}}.
\end{align}
Next, we consider
\[
\mathbf{c}
(\mathbf{u}_h,\mathbf{u}_h;
\mathcal{R}\mathbf{v}_h,\mathcal{R}\mathbf{v}_h).
\]

For every $T\in\mathcal{T}_h$, we first have
\begin{align*}
&
(\mathbf{u}_h\cdot\nabla\mathcal{R}\mathbf{v}_h,
 \mathcal{R}\mathbf{v}_h)_T
\\
&\qquad=
\frac{1}{2}
(\mathbf{u}_h\cdot\nabla
|\mathcal{R}\mathbf{v}_h|^2,1)_T
\\
&\qquad=
\frac{1}{2}
\int_{\partial T}
(\mathbf{u}_h\cdot\mathbf{n}_T)
|\mathcal{R}\mathbf{v}_h|^2
-
\frac{1}{2}
((\nabla\cdot\mathbf{u}_h)\mathcal{R}\mathbf{v}_h,
 \mathcal{R}\mathbf{v}_h)_T.
\end{align*}
Therefore,
\begin{align}
&
(\mathbf{u}_h\cdot\nabla\mathcal{R}\mathbf{v}_h,
 \mathcal{R}\mathbf{v}_h)_T
+
\frac{1}{2}
((\nabla\cdot\mathbf{u}_h)\mathcal{R}\mathbf{v}_h,
 \mathcal{R}\mathbf{v}_h)_T
\nonumber\\
&\qquad=
\frac{1}{2}
\int_{\partial T}
(\mathbf{u}_h\cdot\mathbf{n}_T)
|\mathcal{R}\mathbf{v}_h|^2.
\label{eq:conv_ibp}
\end{align}
Summing over all elements and rearranging the resulting interface
terms gives an equivalent discrete integration-by-parts representation
of the convection form; see
\cite{girault2005splitting,girault2005discontinuous}.
In the present case, this representation reads
\begin{align*}
&
\mathbf{c}
(\mathbf{u}_h,\mathbf{u}_h;
\mathcal{R}\mathbf{v}_h,\mathcal{R}\mathbf{v}_h)
\\
={}&
-
(\mathbf{u}_h\cdot\nabla\mathcal{R}\mathbf{v}_h,
 \mathcal{R}\mathbf{v}_h)_{\mathcal{T}_h}
-
\frac{1}{2}
((\nabla\cdot\mathbf{u}_h)\mathcal{R}\mathbf{v}_h,
 \mathcal{R}\mathbf{v}_h)_{\mathcal{T}_h}
\\
&+
\frac{1}{2}
\left\langle
[\mathbf{u}_h],
\{|\mathcal{R}\mathbf{v}_h|^2\}
\right\rangle_{\mathcal{E}_h^o}
\\
&-
\sum_{T\in\mathcal{T}_h}
\int_{\partial T_-^{\mathbf{u}_h}\setminus\partial\Omega}
|\{\mathbf{u}_h\}\cdot\mathbf{n}_T|
\left(
(\mathcal{R}\mathbf{v}_h)^{\mathrm{int}}
-
(\mathcal{R}\mathbf{v}_h)^{\mathrm{ext}}
\right)
\cdot
(\mathcal{R}\mathbf{v}_h)^{\mathrm{ext}}
\\
&+
\frac{1}{2}
\int_{\partial\Omega}
|\mathbf{u}_h\cdot\mathbf{n}|
|\mathcal{R}\mathbf{v}_h|^2 .
\end{align*}

On the other hand, the original definition gives
\begin{align*}
&
\mathbf{c}
(\mathbf{u}_h,\mathbf{u}_h;
\mathcal{R}\mathbf{v}_h,\mathcal{R}\mathbf{v}_h)
\\
={}&
(\mathbf{u}_h\cdot\nabla\mathcal{R}\mathbf{v}_h,
 \mathcal{R}\mathbf{v}_h)_{\mathcal{T}_h}
+
\frac{1}{2}
((\nabla\cdot\mathbf{u}_h)\mathcal{R}\mathbf{v}_h,
 \mathcal{R}\mathbf{v}_h)_{\mathcal{T}_h}
\\
&-
\frac{1}{2}
\left\langle
[\mathbf{u}_h],
\{|\mathcal{R}\mathbf{v}_h|^2\}
\right\rangle_{\mathcal{E}_h}
\\
&+
\sum_{T\in\mathcal{T}_h}
\int_{\partial T_-^{\mathbf{u}_h}}
|\{\mathbf{u}_h\}\cdot\mathbf{n}_T|
\left(
(\mathcal{R}\mathbf{v}_h)^{\mathrm{int}}
-
(\mathcal{R}\mathbf{v}_h)^{\mathrm{ext}}
\right)
\cdot
(\mathcal{R}\mathbf{v}_h)^{\mathrm{int}} .
\end{align*}

Adding the above two representations, the volume contributions cancel
exactly. Moreover, on the interior edges,
\[
-\frac{1}{2}
\left\langle
[\mathbf{u}_h],
\{|\mathcal{R}\mathbf{v}_h|^2\}
\right\rangle_{\mathcal{E}_h^o}
+
\frac{1}{2}
\left\langle
[\mathbf{u}_h],
\{|\mathcal{R}\mathbf{v}_h|^2\}
\right\rangle_{\mathcal{E}_h^o}
=0.
\]
Thus, the interior part of the consistency term is cancelled exactly
by the corresponding term in the discrete integration-by-parts
representation. The remaining boundary contribution is combined with
the boundary terms according to the boundary convention of the upwind
form.

For the interior upwind contributions, we obtain
\begin{align*}
&
\left(
(\mathcal{R}\mathbf{v}_h)^{\mathrm{int}}
-
(\mathcal{R}\mathbf{v}_h)^{\mathrm{ext}}
\right)
\cdot
(\mathcal{R}\mathbf{v}_h)^{\mathrm{int}}
\\
&\quad-
\left(
(\mathcal{R}\mathbf{v}_h)^{\mathrm{int}}
-
(\mathcal{R}\mathbf{v}_h)^{\mathrm{ext}}
\right)
\cdot
(\mathcal{R}\mathbf{v}_h)^{\mathrm{ext}}
\\
&=
\left|
(\mathcal{R}\mathbf{v}_h)^{\mathrm{int}}
-
(\mathcal{R}\mathbf{v}_h)^{\mathrm{ext}}
\right|^2.
\end{align*}

Consequently,
\begin{align}
&
2\mathbf{c}
(\mathbf{u}_h,\mathbf{u}_h;
\mathcal{R}\mathbf{v}_h,\mathcal{R}\mathbf{v}_h)
\nonumber\\
&=
\sum_{T\in\mathcal{T}_h}
\int_{\partial T_-^{\mathbf{u}_h}}
|\{\mathbf{u}_h\}\cdot\mathbf{n}_T|
\left|
(\mathcal{R}\mathbf{v}_h)^{\mathrm{int}}
-
(\mathcal{R}\mathbf{v}_h)^{\mathrm{ext}}
\right|^2
\nonumber\\
&\quad+
\sum_{T\in\mathcal{T}_h}
\int_{\partial T_-^{\mathbf{u}_h}\cap\partial\Omega}
|\mathbf{u}_h\cdot\mathbf{n}_T|
|\mathcal{R}\mathbf{v}_h|^2 \\
&\geq 0.
\label{eq:conv_upwind_identity}
\end{align}
\end{proof}

\subsection{Error Analysis of the PR-EG Method}

\begin{Lemma}
	$\forall \mathbf{v}_h \in \mathbf{V}_h, q_h \in Q_h$, we have 
	\begin{align}
		\|\mathbf{e}_h\|_{\mathcal{E}} \le Ch (\mu + 1) \| \mathbf{u} \|_2, \\
		\|\epsilon_h\|_{0} \le Ch (\mu + 1) \| \mathbf{u} \|_2.
	\end{align} 
\end{Lemma}
\begin{proof}
	Similar to the error analysis of the standard EG method, we first derive the corresponding error equation.
	Owing to the properties of the reconstruction
	operator $\mathcal{R}$, we have
	\begin{align*}
	(\nabla p,\mathcal{R}\mathbf{v}_h)
	=
	-\mathbf{b}(\mathbf{v}_h,\mathcal{P}_0p).
	\end{align*}
	Hence the pressure approximation error
	$\xi_h=p-\mathcal{P}_0p$ does not enter the error equation.
	\begin{align}
		&\mu \mathbf{a}(\mathbf{e}_h,\mathbf{v}_h)
		-\mathbf{b}(\mathbf{v}_h,\epsilon_h)
		+\mathbf{b}(\mathbf{e}_h,q_h)
		+\mathbf{c}(\mathbf{u},\mathbf{u};\mathbf{u},\mathcal{R} \mathbf{v}_h)
		-\tilde{\mathbf{c}}(\mathbf{u}_h,\mathbf{u}_h;\mathbf{u}_h,\mathbf{v}_h)
		\nonumber\\
		&\qquad=
		-\mu \mathbf{a}(\boldsymbol{\chi}_h,\mathbf{v}_h) + \mu (\Delta \mathbf{u}, \mathcal{R} \mathbf{v}_h - \mathbf{v}_h)
		-\mathbf{b}(\boldsymbol{\chi}_h,q_h).
	\end{align}
	For the term $\mathbf{b}(\boldsymbol{\chi}_h,q_h)$, we have 
	\begin{align*}
		\mathbf{b}(\mathbf{v}-\Pi_h\mathbf{v},q_h)
		=
		0,
		\quad
		\forall q_h\in Q_h.
	\end{align*}
	For the nonlinear terms, we derive
	\begin{align*}
		&\mathbf{c}(\mathbf{u},\mathbf{u};\mathbf{u},\mathcal{R} \mathbf{v}_h)
		-\tilde{\mathbf{c}}(\mathbf{u}_h,\mathbf{u}_h;\mathbf{u}_h,\mathbf{v}_h) \\
		&\qquad=\mathbf{c}(\mathbf{u},\mathbf{u};\mathbf{u},\mathcal{R} \mathbf{v}_h)
		-\mathbf{c}(\mathbf{u}_h,\mathbf{u}_h;\mathcal{R}\mathbf{u}_h,\mathcal{R}\mathbf{v}_h) \\
		&\qquad= \mathbf{c}(\mathbf{u}_h,\mathbf{u};\mathbf{u},\mathcal{R} \mathbf{v}_h)
		-\mathbf{c}(\mathbf{u}_h,\mathbf{u}_h;\mathcal{R}\mathbf{u}_h,\mathcal{R}\mathbf{v}_h) \\
		&\qquad= \mathbf{c}(\mathbf{u}_h,\mathbf{u}_h;\mathcal{R}\mathbf{e}_h,\mathcal{R}\mathbf{v}_h) + \mathbf{c}(\mathbf{u}_h,\mathbf{e}_h;\mathcal{R}\Pi_h\mathbf{u},\mathcal{R}\mathbf{v}_h) + \mathbf{c}(\mathbf{u}_h,\boldsymbol{\chi}_h;\mathcal{R}\Pi_h\mathbf{u},\mathcal{R}\mathbf{v}_h) \nonumber\\
		&\qquad + \mathbf{c}(\mathbf{u}_h,\mathbf{u};\boldsymbol{\chi}_h,\mathcal{R}\mathbf{v}_h) + \mathbf{c}(\mathbf{u}_h,\mathbf{u};\Pi_h\mathbf{u}-\mathcal{R}\Pi_h\mathbf{u},\mathcal{R}\mathbf{v}_h).
	\end{align*}

	Therefore, from the inf-sup condition \eqref{eq:inf-sup2}, we can obtain
	\begin{align}
		C_{IS}
		\|(\mathbf{e}_h,\epsilon_h)\|_{\mathbf{V}_h\times Q_h}
		&\leq
		\sup_{(\mathbf{v}_h,q_h)\in\mathbf{V}_h\times Q_h}
		\frac{
		\mathcal{S}_{\mathbf{u}_h}
		((\mathbf{e}_h,\epsilon_h),(\mathbf{v}_h,q_h))
		}{
		\|(\mathbf{v}_h,q_h)\|_{\mathbf{V}_h\times Q_h}
		} \\ \nonumber
		&\le \sup_{(\mathbf{v}_h,q_h)\in\mathbf{V}_h\times Q_h} \frac{1}{\|(\mathbf{v}_h,q_h)\|_{\mathbf{V}_h\times Q_h}} ( -\mu \mathbf{a}(\boldsymbol{\chi}_h,\mathbf{v}_h) + \mu (\Delta \mathbf{u}, \mathcal{R} \mathbf{v}_h - \mathbf{v}_h) \\ 
		&- \mathbf{c}(\mathbf{u}_h,\mathbf{e}_h;\mathcal{R}\Pi_h\mathbf{u},\mathcal{R}\mathbf{v}_h) - \mathbf{c}(\mathbf{u}_h,\boldsymbol{\chi}_h;\mathcal{R}\Pi_h\mathbf{u},\mathcal{R}\mathbf{v}_h) \nonumber \\
		&- \mathbf{c}(\mathbf{u}_h,\mathbf{u};\boldsymbol{\chi}_h,\mathcal{R}\mathbf{v}_h) - \mathbf{c}(\mathbf{u}_h,\mathbf{u};\Pi_h\mathbf{u}-\mathcal{R}\Pi_h\mathbf{u},\mathcal{R}\mathbf{v}_h)). \label{PREG_inf_sup_error}
	\end{align}
	The first term in \eqref{PREG_inf_sup_error} can be estimated as follows:
	\begin{align*}
		\mu \mathbf{a}(\boldsymbol{\chi}_h,\mathbf{v}_h) \le \mu Ch \| \mathbf{u} \|_2 \| \mathbf{v}_h \|_{\mathcal{E}}.
	\end{align*}
	Besides, the Cauchy-Schwarz inequality and the second term lead to 
	\begin{align}
		\mu (\Delta \mathbf{u}, \mathcal{R} \mathbf{v}_h - \mathbf{v}_h) &\le \mu \| \mathbf{u} \|_2 \| \mathcal{R} \mathbf{v}_h - \mathbf{v}_h \|_0 \\
		& \le \mu Ch \| \mathbf{u} \|_2 \| \mathbf{v}_h \|_{\mathcal{E}}.
	\end{align}
	The nonlinear terms can be estimated as follows:
	\begin{align}
		& \mathbf{c}(\mathbf{u}_h,\boldsymbol{\chi}_h;\mathcal{R}\Pi_h\mathbf{u},\mathcal{R}\mathbf{v}_h) + \mathbf{c}(\mathbf{u}_h,\mathbf{u};\boldsymbol{\chi}_h,\mathcal{R}\mathbf{v}_h) + \mathbf{c}(\mathbf{u}_h,\mathbf{u};\Pi_h\mathbf{u}-\mathcal{R}\Pi_h\mathbf{u},\mathcal{R}\mathbf{v}_h) \nonumber \\ 
		&\qquad \le C_0^{PR} ( \| \boldsymbol{\chi}_h \|_{\mathcal{E}} \| \mathcal{R}\Pi_h\mathbf{u} \|_{\mathcal{E}} \| \mathcal{R}\mathbf{v}_h \|_{\mathcal{E}} + \| \mathbf{u} \|_{\mathcal{E}} \| \boldsymbol{\chi}_h \|_{\mathcal{E}} \| \mathcal{R}\mathbf{v}_h \|_{\mathcal{E}} + \| \mathbf{u} \|_{\mathcal{E}} \| \Pi_h\mathbf{u}-\mathcal{R}\Pi_h\mathbf{u} \|_{\mathcal{E}} \| \mathcal{R}\mathbf{v}_h \|_{\mathcal{E}} ) \nonumber \\
		&\qquad \le C_0^{PR} ( Ch \| \mathbf{u} \|_2 | \mathbf{u} |_1 + Ch \| \mathbf{u} \|_2 + C\| \Pi_h\mathbf{u}-\mathcal{R}\Pi_h\mathbf{u} \|_{\mathcal{E}})\| \mathbf{v}_h \|_{\mathcal{E}}. \label{nonlinear_error}
	\end{align}
	To complete the estimate, it remains to bound the reconstruction error term
	$\|\Pi_h\mathbf{u}-\mathcal{R}\Pi_h\mathbf{u}\|_{\mathcal{E}}$.
	From the definition of the $\mathcal{E}$-norm, we have
	\begin{align*}
		\|\Pi_h\mathbf{u}-\mathcal{R}\Pi_h\mathbf{u}\|_{\mathcal{E}}^2
		&=
		\|\nabla(\Pi_h\mathbf{u}-\mathcal{R}\Pi_h\mathbf{u})
		\|_{0,\mathcal{T}_h}^2 + \rho\|h_e^{-1/2}[\Pi_h\mathbf{u}-\mathcal{R}\Pi_h\mathbf{u}]\|_{0,\mathcal{E}_h}^2.
	\end{align*}

	For the first term, using the inverse inequality together with \eqref{eq:scaling argument}, we obtain
	\begin{align*}
		\|\nabla(\Pi_h\mathbf{u}-\mathcal{R}\Pi_h\mathbf{u})
		\|_{0,\mathcal{T}_h}^2
		&\leq
		C \sum_{K \in \mathcal{T}_h} h_{K}^{-2} \|\Pi_h\mathbf{u}-\mathcal{R}\Pi_h\mathbf{u}\|_{0,K}^2 \\
		&\leq
		C \sum_{e \in \mathcal{E}_h} h_{K}^{-2} h_e \|[\Pi_h\mathbf{u}]\|_{0,e}^2 \\
		&\leq
		C
		\|h_e^{-1/2}[\Pi_h\mathbf{u}]
		\|_{0,\mathcal{E}_h}^2.
	\end{align*}
	Since $\mathbf{u}\in[H_0^1(\Omega)]^d$, we have
	$[\mathbf{u}]=0$ on $\mathcal{E}_h$. Hence,
	\begin{align*}
		\|h_e^{-1/2}[\Pi_h\mathbf{u}]
		\|_{0,\mathcal{E}_h}
		&=
		\|h_e^{-1/2}[\Pi_h\mathbf{u}-\mathbf{u}]
		\|_{0,\mathcal{E}_h}
		\\
		&\leq
		C\|\mathbf{u}-\Pi_h\mathbf{u}\|_{\mathcal{E}}
		\\
		&\leq
		Ch\|\mathbf{u}\|_2.
	\end{align*}
	Therefore,
	\begin{align*}
		\|\nabla(\Pi_h\mathbf{u}-\mathcal{R}\Pi_h\mathbf{u})
		\|_{0,\mathcal{T}_h}^2
		\leq
		Ch^2\|\mathbf{u}\|_2^2.
	\end{align*}

	For the second term, the inverse trace inequality and \eqref{eq:scaling argument} give
	\begin{align*}
		\|h_e^{-1/2}
		[\Pi_h\mathbf{u}-\mathcal{R}\Pi_h\mathbf{u}]
		\|_{0,\mathcal{E}_h}
		&\leq
		\|h_e^{-1}
		[\Pi_h\mathbf{u}-\mathcal{R}\Pi_h\mathbf{u}]
		\|_{0,\mathcal{E}_h} \\
		&\leq
		C
		\|h_e^{-1/2}[\Pi_h\mathbf{u}]
		\|_{0,\mathcal{E}_h}
		\\
		&\leq
		Ch\|\mathbf{u}\|_2.
	\end{align*}

	Combining the above estimates, we obtain
	\begin{align*}
		\|\Pi_h\mathbf{u}-\mathcal{R}\Pi_h\mathbf{u}\|_{\mathcal{E}}^2
		\leq
		Ch^2\|\mathbf{u}\|_2^2.
	\end{align*}

	Therefore, \eqref{nonlinear_error} leads to
	\begin{align*}
		& \mathbf{c}(\mathbf{u}_h,\boldsymbol{\chi}_h;\mathcal{R}\Pi_h\mathbf{u},\mathcal{R}\mathbf{v}_h) + \mathbf{c}(\mathbf{u}_h,\mathbf{u};\boldsymbol{\chi}_h,\mathcal{R}\mathbf{v}_h) + \mathbf{c}(\mathbf{u}_h,\mathbf{u};\Pi_h\mathbf{u}-\mathcal{R}\Pi_h\mathbf{u},\mathcal{R}\mathbf{v}_h) \nonumber \\ 
		&\qquad \le C_0^{PR} ( Ch \| \mathbf{u} \|_2 | \mathbf{u} |_1 + Ch \| \mathbf{u} \|_2 + Ch\| \mathbf{u} \|_2)\| \mathbf{v}_h \|_{\mathcal{E}} \\
		&\qquad \le Ch \| \mathbf{u} \|_2 \| \mathbf{v}_h \|_{\mathcal{E}}.
	\end{align*}

	Finally, the remaining nonlinear term can be estimated as follows:
	\begin{align*}
		\mathbf{c}(\mathbf{u}_h,\mathbf{e}_h;\mathcal{R}\Pi_h\mathbf{u},\mathcal{R}\mathbf{v}_h) \le C_0^{PR} \| \mathbf{e}_h \|_{\mathcal{E}} \| \mathcal{R}\Pi_h\mathbf{u} \|_{\mathcal{E}} \| \mathbf{v}_h \|_{\mathcal{E}} \le C \| \mathbf{e}_h \|_{\mathcal{E}} \| \mathbf{v}_h \|_{\mathcal{E}}.
	\end{align*}

	In conclusion, we have
	\begin{align*}
		\|(\mathbf{e}_h,\epsilon_h)\|_{\mathbf{V}_h\times Q_h} \le Ch (\mu + 1)\| \mathbf{u} \|_2 + C \| \mathbf{e}_h \|_{\mathcal{E}}.
	\end{align*}

	Thus, the velocity and pressure error estimates of PR-EG method can be obtained separately as
	\begin{align*}
		\|\mathbf{e}_h\|_{\mathcal{E}} \le Ch (\mu + 1) \| \mathbf{u} \|_2, \\
		\|\epsilon_h\|_{0} \le Ch (\mu + 1) \| \mathbf{u} \|_2.
	\end{align*}
\end{proof}
\begin{Theorem}
	We assume that $(\mathbf{u}, p) \in \textbf{H}^1_0(\Omega) \times L^2_0(\Omega)$ is the solution of Navier-Stokes equations and $(\mathbf{u}_h, p_h)$ is the solution of \eqref{eq:variation euqation_PREG}, we have the following error estimates:
	\begin{align}
		&\| \mathbf{u}-\mathbf{u}_h \|_{\mathcal{E}}
		\leq
		Ch (\mu + 1) \| \mathbf{u} \|_2, \\
		&\| p-p_h \|_0
		\leq
		Ch
		\left(
		(\mu + 1) \| \mathbf{u} \|_2
		+
		\| p \|_1
		\right).
	\end{align}
\end{Theorem}
\begin{proof}
	We first decompose the total error into the interpolation error and the discrete error. By the triangle inequality, the velocity error satisfies
	\begin{align*}
		\| \mathbf{u}-\mathbf{u}_h \|_{\mathcal{E}}
		\leq
		\| \boldsymbol{\chi}_h \|_{\mathcal{E}}
		+
		\| \mathbf{e}_h \|_{\mathcal{E}}.
	\end{align*}

	For the interpolation error, the approximation property of $\Pi_h$ gives
	\begin{align*}
		\| \boldsymbol{\chi}_h \|_{\mathcal{E}}
		\leq
		Ch\| \mathbf{u} \|_2.
	\end{align*}
	Moreover, the discrete error estimate yields
	\begin{align*}
		\| \mathbf{e}_h \|_{\mathcal{E}}
		\leq
		Ch (\mu + 1) \| \mathbf{u} \|_2.
	\end{align*}
	Combining the above estimates, we obtain
	\begin{align}
		\| \mathbf{u}-\mathbf{u}_h \|_{\mathcal{E}}
		\leq
		Ch (\mu + 1) \| \mathbf{u} \|_2.
	\end{align}

	Similarly, we can obtain the error estimate for pressure in the same way
	\begin{align}
		\| p-p_h \|_0
		\leq
		Ch
		\left(
		(\mu + 1) \| \mathbf{u} \|_2
		+
		\| p \|_1
		\right).
	\end{align}
\end{proof}

\section{Numerical Experiments}

In this section, we present numerical experiments to illustrate the performance of the Enriched Galerkin (EG) method and PR-EG method for the stationary Navier-Stokes equations. 
All computations are carried out on shape-regular triangular meshes obtained by uniform refinement of an initial partition of the computational domain. 
The EG velocity space $V_h$ and pressure space $Q_h$ are defined in Section~3, and we denote by $(u_h,p_h)\in V_h\times Q_h$ the discrete solution.

\subsection{Numerical Experiments for EG method}
Let $h := \max_{T\in\mathcal{T}_h} h_T$ be the mesh size. 
The penalty parameter $\rho$ is fixed to $\rho = 10$, which is sufficiently large to guarantee the coercivity of the bilinear form associated with the symmetric interior penalty discretization of the viscous term. 
The viscosity is set to $\mu=1$.
On each mesh level, the nonlinear EG system is solved by a Picard iteration, starting from the initial guess $(u_h^{(0)},p_h^{(0)})=(0,0)$.

\subsubsection{Example 1: Smooth Manufactured Solution for EG method}

In the first example we consider a smooth manufactured solution on the unit square $\Omega=(0,1)^2$. The exact velocity and pressure are given by
\begin{align*}
	u_1(x,y) &= 2 x^2 (1-x)^2\, y(1-y)\,(1-2y),\\
	u_2(x,y) &= -2 y^2 (1-y)^2\, x(1-x)\,(1-2x),\\
	p(x,y)   &= \sin(\pi x)\cos(\pi y).
\end{align*}
A direct calculation shows that $\nabla\cdot u = 0$ in $\Omega$ and that the pressure has zero mean over $\Omega$, so that $p\in L^2_0(\Omega)$.\\
\indent The initial mesh is a uniform triangulation of $\Omega$ into right triangles, and subsequent meshes are obtained by uniform refinement, so that $h$ is reduced by a factor of $1/2$ on each level. 
On each mesh, we compute the EG approximation $(u_h,p_h)$ of problem \eqref{eq:model_problem}.\\
\indent Since the exact velocity and pressure are smooth polynomials (with $u\in H^2(\Omega)^2$
and $p\in H^1(\Omega)$), the interpolation estimates for the EG space and the a priori
error bounds derived in Section~4 suggest the following asymptotic behavior as $h\to0$:
\[
\| u - u_h \|_{\mathcal{E}}
\;\lesssim\; Ch
		\left(
		(\mu+1)\| \mathbf{u} \|_2
		+
		\| p \|_1
		\right), \qquad
\|p-p_h\|_{0}
\;\lesssim\; Ch
		\left(
		(\mu+1)\| \mathbf{u} \|_2
		+
		\| p \|_1
		\right).
\]

With a constant $C$ independent of $h$, the numerical results reported in Table~\ref{tab:EG_poly_errors} confirm the theoretical convergence rates. \\
In particular, the velocity approximation achieves first-order accuracy in the mesh-dependent energy norm $\| u - u_h \|_{\mathcal{E}}$ and second-order accuracy in the $L^2$-norm $\|u - u_h\|_{0}$. \\
The pressure approximation converges with first-order accuracy in the $L^2$-norm $\|p - p_h\|_{0}$, which is fully consistent with the regularity of the exact solution and the analysis presented in Section~4. \\

\begin{table}[H]
	\centering
	\caption{A mesh refinement study for the velocity errors of EG with h when $\mu = 1$.}
	\label{tab:EG_poly_errors}
	\begin{tabular}{c|cc|cc|cc}
		\hline
		$h$ 
		& $\| u - u_h \|_{0}$ & Order 
		& $\|u-u_h\|_{\mathcal{E}}$ & Order
		& $\|p-p_h\|_{0}$ & Order \\
		\hline
		$1/4$  & 1.6089e$-$01 &   --      & 4.1265e$-$01 &   --      & 1.9456e$+$00 &   --      \\
		$1/8$  & 3.9909e$-$02 & 2.0113    & 2.2501e$-$01 & 0.8750    & 1.0105e$+$00 & 0.9451    \\
		$1/16$ & 1.0065e$-$02 & 1.9874    & 1.0826e$-$01 & 1.0555    & 5.1951e$-$01 & 0.9598    \\
		$1/32$ & 2.5211e$-$03 & 1.9972    & 5.0424e$-$02 & 1.1023    & 2.6007e$-$01 & 0.9983    \\
		$1/64$ & 6.3008e$-$04 & 2.0004    & 2.4120e$-$02 & 1.0639    & 1.3310e$-$01 & 0.9664    \\
		\hline
	\end{tabular}
\end{table}

\begin{figure}[h]
	\centering
	\begin{subcaptiongroup}
		\begin{subfigure}{0.45\textwidth}
			\centering
			\includegraphics[width=\textwidth]{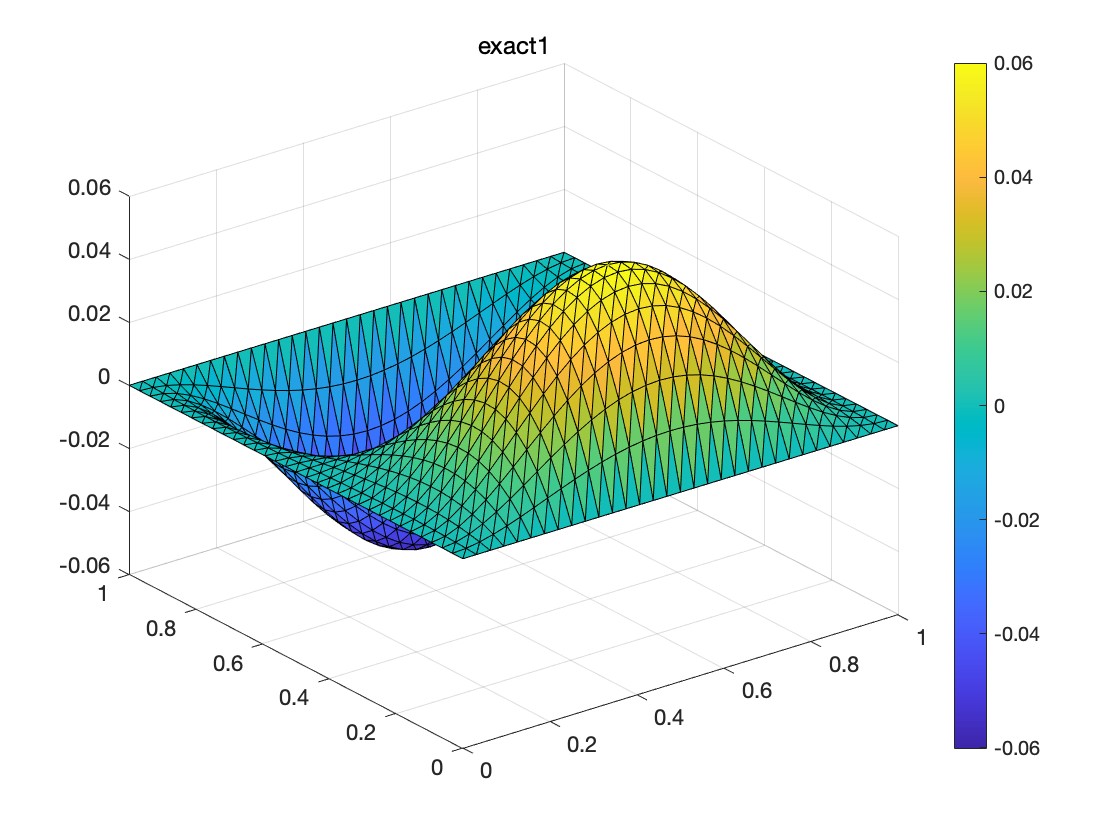} 
			\caption{Exact velocity component $u_1$.}
			\label{fig:EG_poly_u1_exact}
		\end{subfigure}
		\hfill
		\begin{subfigure}{0.45\textwidth}
			\centering
			\includegraphics[width=\textwidth]{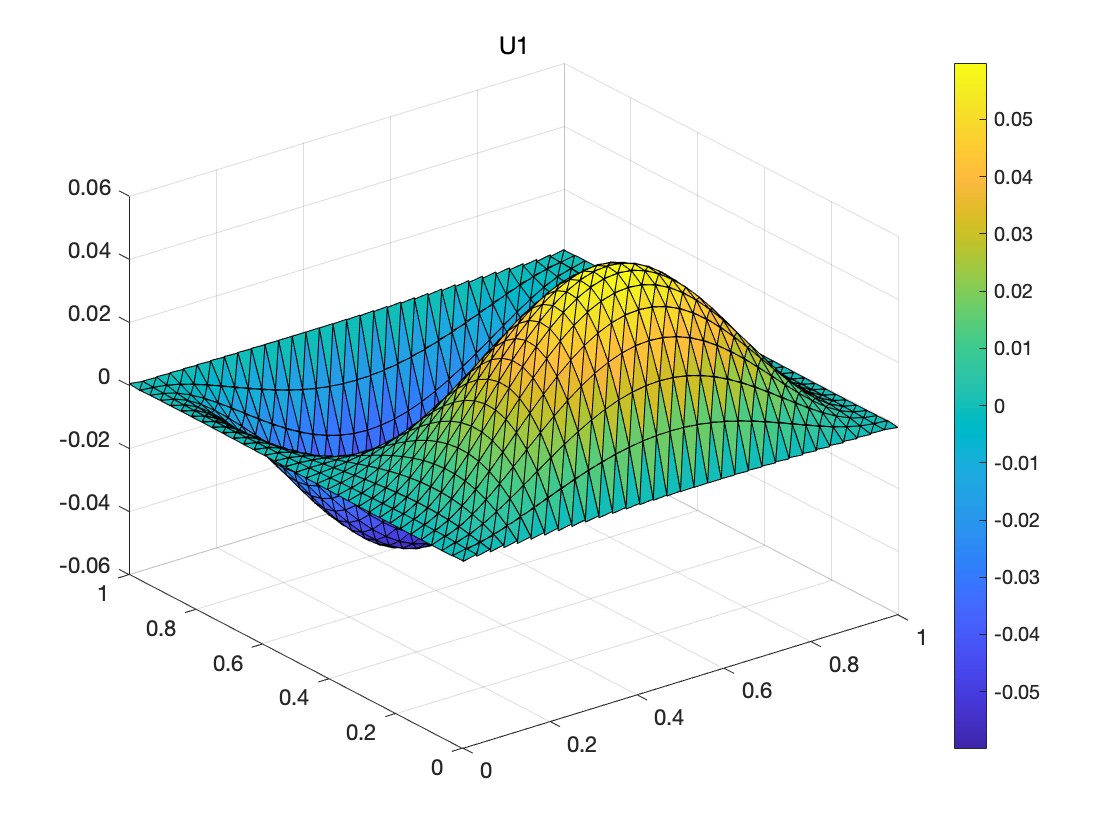} 
			\caption{EG approximation $u_{1,h}$.}
			\label{fig:EG_poly_u1_num}
		\end{subfigure}
		
		\vskip\baselineskip
		
		\begin{subfigure}{0.45\textwidth}
			\centering
			\includegraphics[width=\textwidth]{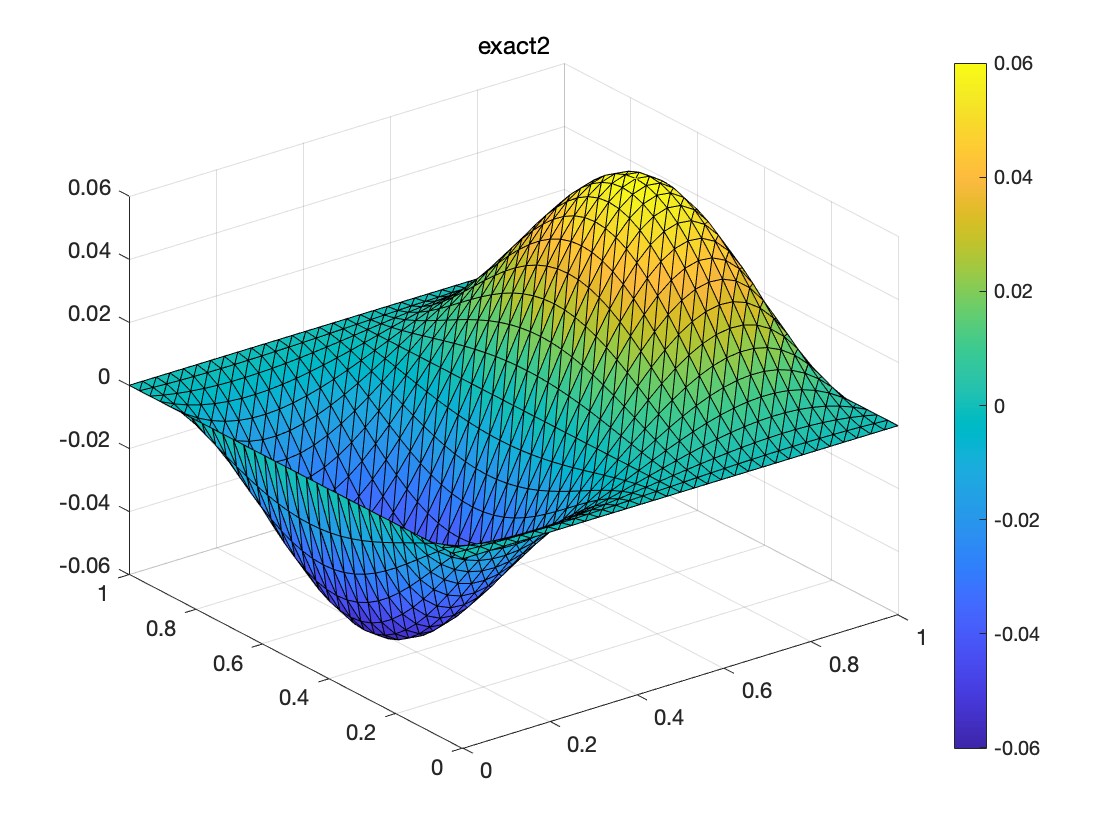} 
			\caption{Exact velocity component $u_2$.}
			\label{fig:EG_poly_u2_exact}
		\end{subfigure}
		\hfill
		\begin{subfigure}{0.45\textwidth}
			\centering
			\includegraphics[width=\textwidth]{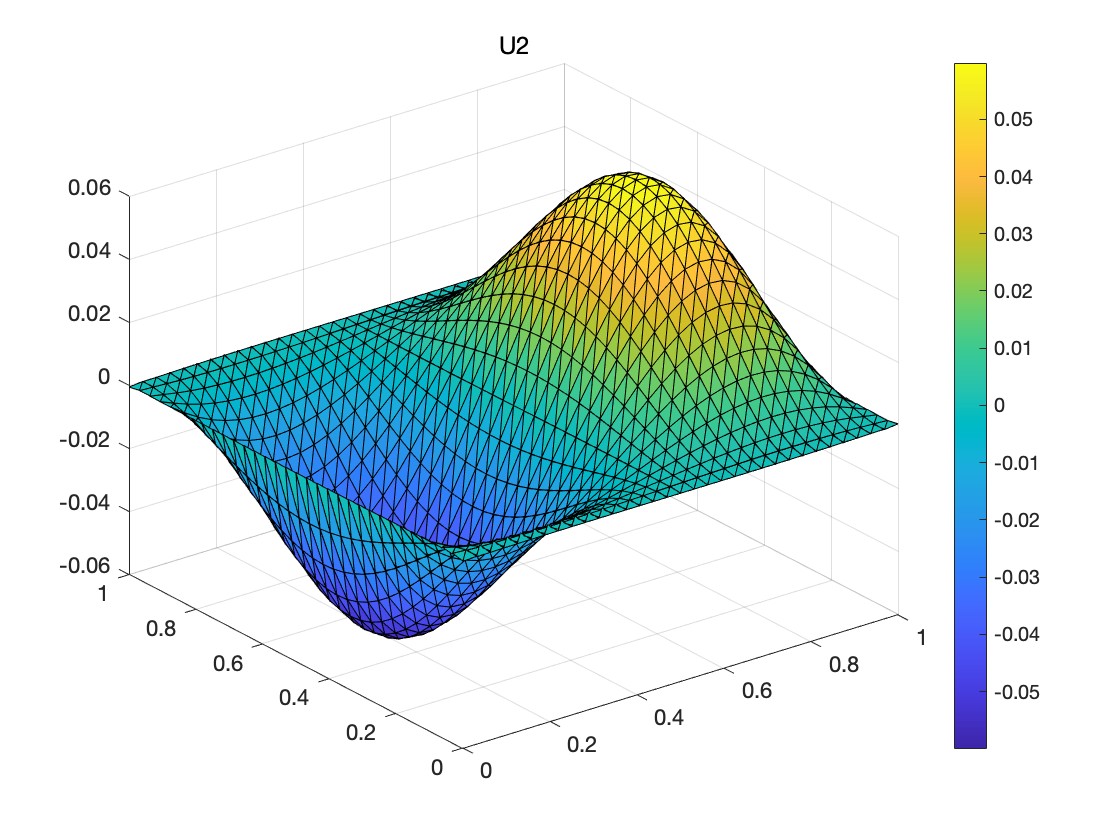} 
			\caption{EG approximation $u_{2,h}$.}
			\label{fig:EG_poly_u2_num}
		\end{subfigure}
	\end{subcaptiongroup}
	\caption{Exact and EG velocity components for Example~1 (smooth polynomial solution).}
	\label{fig:EG_poly_u}
\end{figure}

This experiment demonstrates that, for a smooth incompressible flow, the enriched
Galerkin discretization of the stationary Navier--Stokes equations achieves the
expected optimal convergence orders in the mesh-dependent energy norm and in the
standard $L^2$-norms for both velocity and pressure.

\subsubsection{Example 2: Lid-driven Cavity Flow for EG method}

\begin{figure}[H] 
	\begin{flushleft}
		\includegraphics[width=1.0\textwidth]{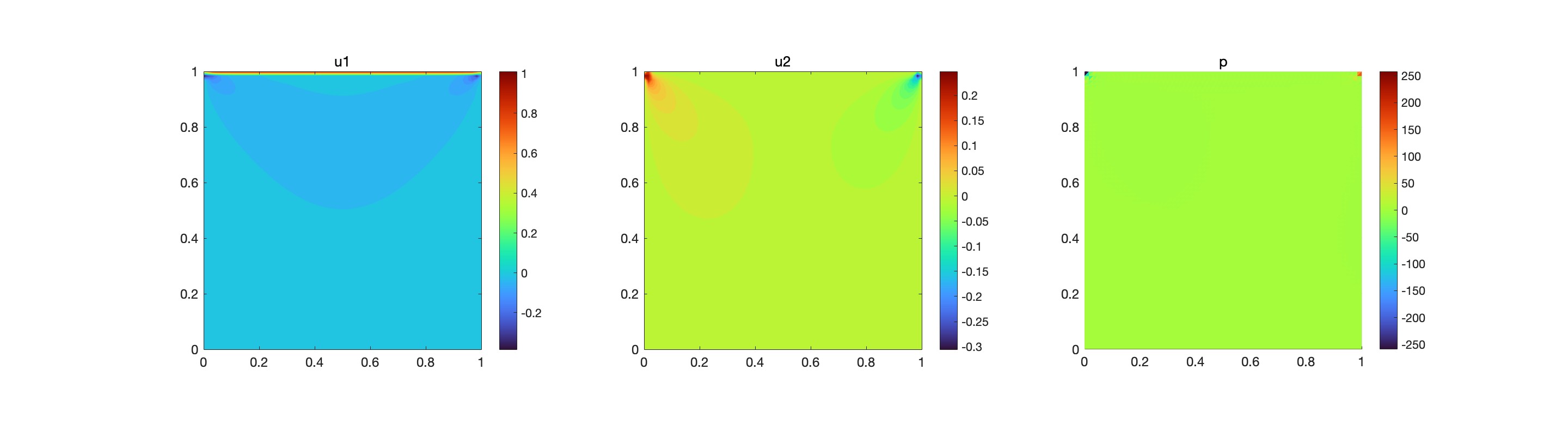} 
		\caption{EG approximation of the lid-driven cavity flow: contour plots of $u_{1,h}$, $u_{2,h}$, and $p_h$ on a mesh with $h_0 = 1/32$.} 
		\label{fig:cavity_plot} 
	 \end{flushleft}
\end{figure}

In this example, we consider the classical lid-driven cavity flow in the unit square $\Omega = (0,1)^2$. 
The body force is set to zero. 
On the top boundary, we prescribe a tangential unit velocity $u = (1,0)^{T}$, while homogeneous Dirichlet boundary conditions are imposed on the remaining three sides. 
The pressure is determined up to a constant and is fixed by enforcing a zero-mean condition. The EG formulation is assembled exactly as in the previous example. 
To obtain an initial guess for the nonlinear iteration, we first solve the Stokes problem by dropping the convective term, and then the full Navier--Stokes equations are solved using a Picard iteration. 
The iteration stops once the relative $L^2$ difference between successive iterates satisfies \[ \frac{\|x^{(k+1)} - x^{(k)}\|_2}{\|x^{(k)}\|_2} < 10^{-10}, \] or a maximum of 20 iterations is reached. After obtaining the discrete velocity and pressure, the elementwise bubble contributions are added back to the nodal values and the solution is interpolated onto a structured Cartesian grid. 
Figure~\ref{fig:cavity_plot} displays the contour plots of the horizontal velocity $u_{1,h}$, the vertical velocity $u_{2,h}$, and the pressure $p_h$ on the mesh with $h_0 = 1/32$. 
Characteristic velocity and pressure structures are captured.

\subsection{Numerical Experiments for the PR-EG Method}

In this subsection, we investigate the accuracy and pressure robustness of the proposed PR-EG method for the stationary Navier--Stokes equations. 
In particular, we compare the standard EG and PR-EG methods for different values of the viscosity parameter $\mu$. The same finite element spaces, mesh family, penalty parameter, and nonlinear solver are employed for both methods so that the effect of the velocity reconstruction can be examined directly.\\
\indent Unless otherwise specified, the penalty parameter is fixed as $\rho=10$.
The nonlinear systems are solved by Newton's method. On each mesh level, the iteration is terminated when the prescribed stopping criterion is satisfied.

\subsubsection{Example 3: Accuracy and Robustness Test with $\mu>0$}

\begin{figure}[htbp]
    \centering

    \begin{subfigure}{0.47\textwidth}
        \centering

        %
        \includegraphics[width=\textwidth]{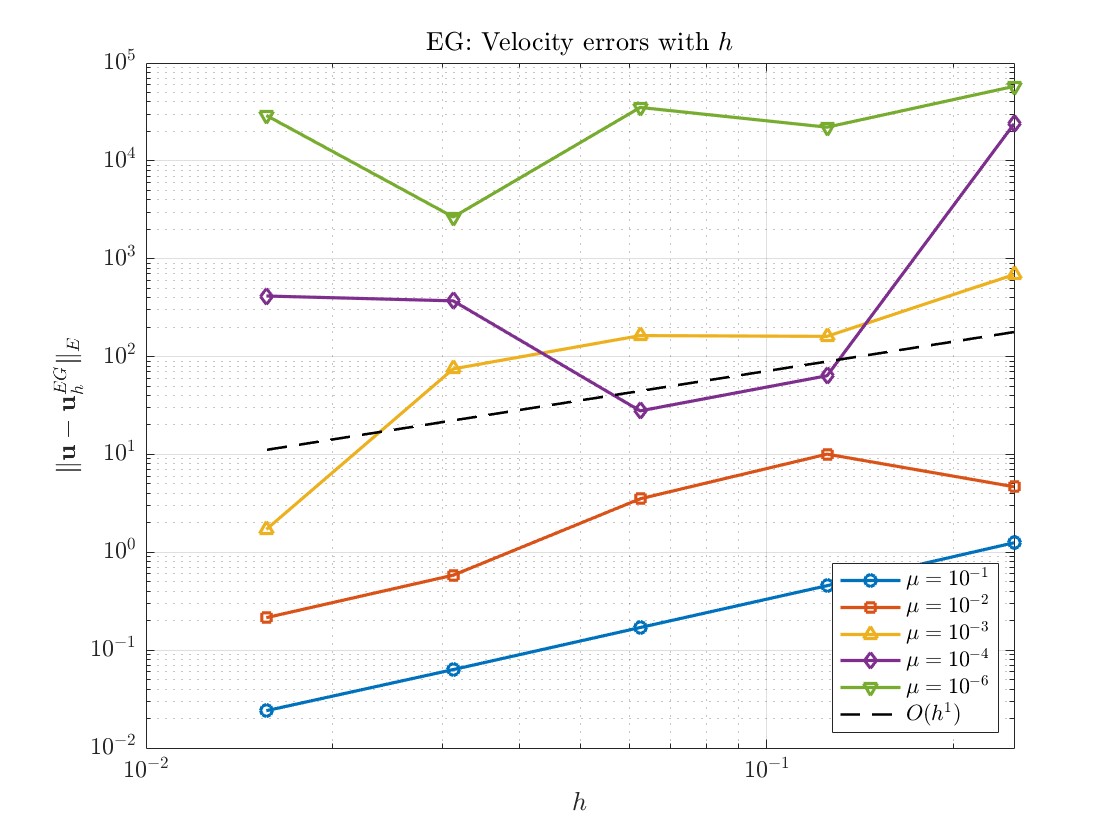}

        \caption{EG: velocity errors for different values of $\mu$.}
        \label{fig:EG_u_error}
    \end{subfigure}
    \hfill
    \begin{subfigure}{0.47\textwidth}
        \centering

        %
        \includegraphics[width=\textwidth]{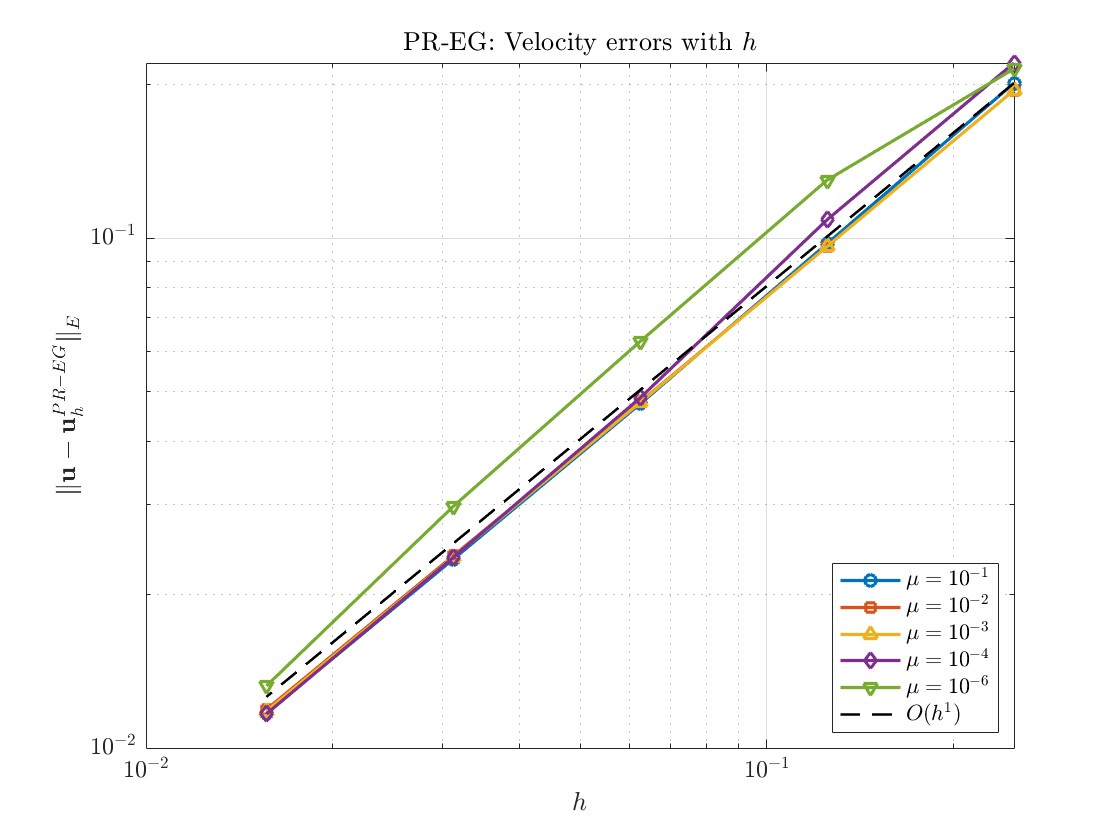}

        \caption{PR-EG: velocity errors for different values of $\mu$.}
        \label{fig:PREG_u_error}
    \end{subfigure}

    \caption{Convergence histories of the velocity errors with respect to
    $h$ for different values of the viscosity parameter $\mu$.}
    \label{fig:preg_mu_h_errors}
\end{figure}

We first consider the smooth manufactured solution on the unit square
\[
    \Omega=(0,1)^2.
\]
The exact velocity and pressure are given by
\begin{align*}
    u_1(x,y)
    &=
    2x^2(1-x)^2y(1-y)(1-2y),\\
    u_2(x,y)
    &=
    -2y^2(1-y)^2x(1-x)(1-2x),\\
    p(x,y)
    &=
    \sin(\pi x)\cos(\pi y).
\end{align*}
A direct calculation shows that
\[
    \nabla\cdot\mathbf u=0
    \qquad\text{in }\Omega,
\]
and that the pressure has zero mean over $\Omega$. The forcing term is
determined by substituting the exact solution into the stationary
Navier--Stokes equations for each value of $\mu$, and the boundary condition
is prescribed according to the exact velocity.\\
\indent The initial mesh is a uniform triangulation of $\Omega$ into right triangles,
and the subsequent meshes are generated by uniform refinement. Thus, the mesh
size $h$ is reduced by a factor of $1/2$ on each refinement level. On every
mesh, both the standard EG approximation
$(\mathbf u_h^{\mathrm{EG}},p_h^{\mathrm{EG}})$ and the PR-EG approximation
$(\mathbf u_h^{\mathrm{PR}},p_h^{\mathrm{PR}})$ are computed.\\
\indent The purpose of this example is twofold. First, we verify the convergence
behavior of the PR-EG method under mesh refinement. Second, we investigate
the sensitivity of the velocity approximation to the viscosity parameter and
compare the behavior of the standard EG and PR-EG methods.\\
\indent For the standard EG method, the velocity error estimate obtained in
Section~4 contains a contribution from the pressure approximation,
\[
    \|\mathbf u-\mathbf u_h^{\mathrm{EG}}\|_{\mathcal E}
    \lesssim
    Ch\left((\mu+1)\|\mathbf u\|_2+\|p\|_1\right).
\]
In contrast, the pressure-dependent contribution is eliminated from the
velocity estimate of the PR-EG method. According to Theorem~5.1,
\[
    \|\mathbf u-\mathbf u_h^{\mathrm{PR}}\|_{\mathcal E}
    \lesssim
    Ch(\mu+1)\|\mathbf u\|_2.
\]
Hence, the reconstructed formulation is expected to provide a velocity
approximation that is less sensitive to pressure effects, particularly when
the viscous contribution becomes small.\\
\indent To examine this behavior numerically, we take
\[
    \mu=10^{-1},\;10^{-2},\;10^{-3},\;10^{-4},\;10^{-6},
\]
and compute the velocity errors on a sequence of uniformly refined meshes.
Figure~\ref{fig:preg_mu_h_errors} compares the convergence histories of the
standard EG and PR-EG methods. The horizontal axis denotes the mesh size $h$,
while the vertical axis represents the velocity error. A reference line with
slope one is included to indicate the theoretically predicted first-order
behavior in the mesh-dependent energy norm.

To provide a more detailed comparison in the small-viscosity regime, we fix
\[
    \mu=10^{-6}
\]
and report the numerical errors and convergence rates for both methods in
Table~\ref{tab:eg_preg_mu_small}. We measure the velocity error in the
mesh-dependent energy norm and the $L^2$-norm, together with the pressure
error in the $L^2$-norm.

\begin{table}[htbp]
\centering
\caption{A mesh refinement study for the EG and PR-EG methods when
$\mu=10^{-6}$.}
\label{tab:eg_preg_mu_small}

\begin{tabular}{c|cc|cc|cc}
    \hline
    \multicolumn{7}{c}{Standard EG method}\\
    \hline
    $h$
    & $\|\mathbf u-\mathbf u_h^{\mathrm{EG}}\|_{\mathcal E}$
    & Order
    & $\|\mathbf u-\mathbf u_h^{\mathrm{EG}}\|_0$
    & Order
    & $\|p-p_h^{\mathrm{EG}}\|_0$
    & Order\\
    \hline
    $1/4$
    & $2.3297\mathrm{e}{+04}$ & --
    & $5.0304\mathrm{e}{+03}$ & --
    & $1.5479\mathrm{e}{+00}$ & -- \\

    $1/8$
    & $5.3192\mathrm{e}{+03}$ & $2.1308$
    & $5.2220\mathrm{e}{+02}$ & $3.2680$
    & $2.1758\mathrm{e}{+01}$ & $-3.8131$ \\

    $1/16$
    & $5.2184\mathrm{e}{+02}$ & $3.3495$
    & $2.7658\mathrm{e}{+01}$ & $4.2388$
    & $1.0973\mathrm{e}{+01}$ & $0.9876$ \\

    $1/32$
    & $3.9166\mathrm{e}{+03}$ & $-2.9079$
    & $9.9005\mathrm{e}{+01}$ & $-1.8398$
    & $1.1788\mathrm{e}{+02}$ & $-3.4254$ \\

    $1/64$
    & $1.2245\mathrm{e}{+04}$ & $-1.6445$
    & $1.5893\mathrm{e}{+02}$ & $-0.6828$
    & $7.4730\mathrm{e}{+01}$ & $0.6576$ \\
    \hline

    \multicolumn{7}{c}{}\\[-2mm]

    \hline
    \multicolumn{7}{c}{PR-EG method}\\
    \hline
    $h$
    & $\|\mathbf u-\mathbf u_h^{\mathrm{PR}}\|_{\mathcal E}$
    & Order
    & $\|\mathbf u-\mathbf u_h^{\mathrm{PR}}\|_0$
    & Order
    & $\|p-p_h^{\mathrm{PR}}\|_0$
    & Order\\
    \hline
    $1/4$
    & $2.1494\mathrm{e}{-01}$ & --
    & $2.4831\mathrm{e}{-02}$ & --
    & $9.5470\mathrm{e}{-01}$ & -- \\

    $1/8$
    & $1.3001\mathrm{e}{-01}$ & $0.7253$
    & $1.4388\mathrm{e}{-02}$ & $0.7873$
    & $4.8019\mathrm{e}{-01}$ & $0.9915$ \\

    $1/16$
    & $6.2730\mathrm{e}{-02}$ & $1.0514$
    & $3.6585\mathrm{e}{-03}$ & $1.9755$
    & $2.4045\mathrm{e}{-01}$ & $0.9979$ \\

    $1/32$
    & $2.9779\mathrm{e}{-02}$ & $1.0749$
    & $9.2378\mathrm{e}{-04}$ & $1.9856$
    & $1.2027\mathrm{e}{-01}$ & $0.9995$ \\

    $1/64$
    & $1.3235\mathrm{e}{-02}$ & $1.1700$
    & $2.3024\mathrm{e}{-04}$ & $2.0044$
    & $6.0139\mathrm{e}{-02}$ & $0.9999$ \\
    \hline
\end{tabular}

\end{table}

\begin{figure}[htbp]
    \centering

    \begin{subfigure}{0.45\textwidth}
        \centering
        \includegraphics[width=\textwidth]{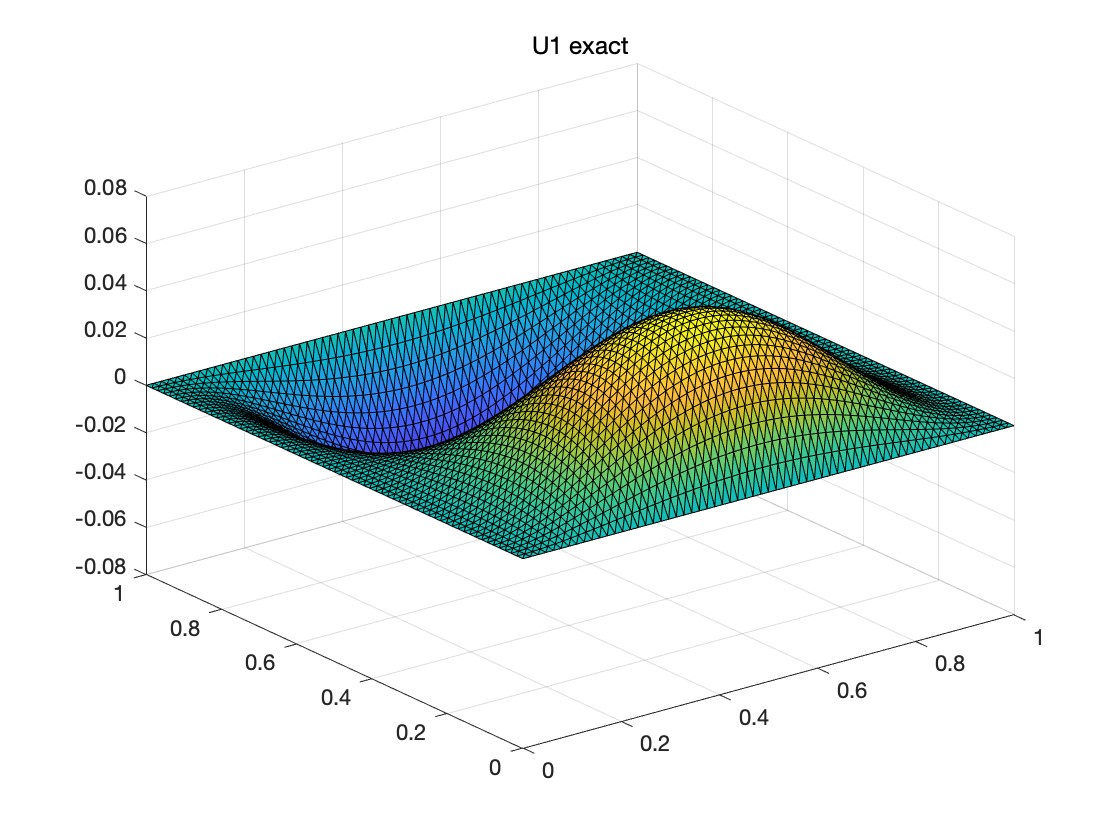}
        \caption{Exact velocity component $u_1$.}
        \label{fig:preg_u1_exact}
    \end{subfigure}
    \hfill
    \begin{subfigure}{0.45\textwidth}
        \centering
        \includegraphics[width=\textwidth]{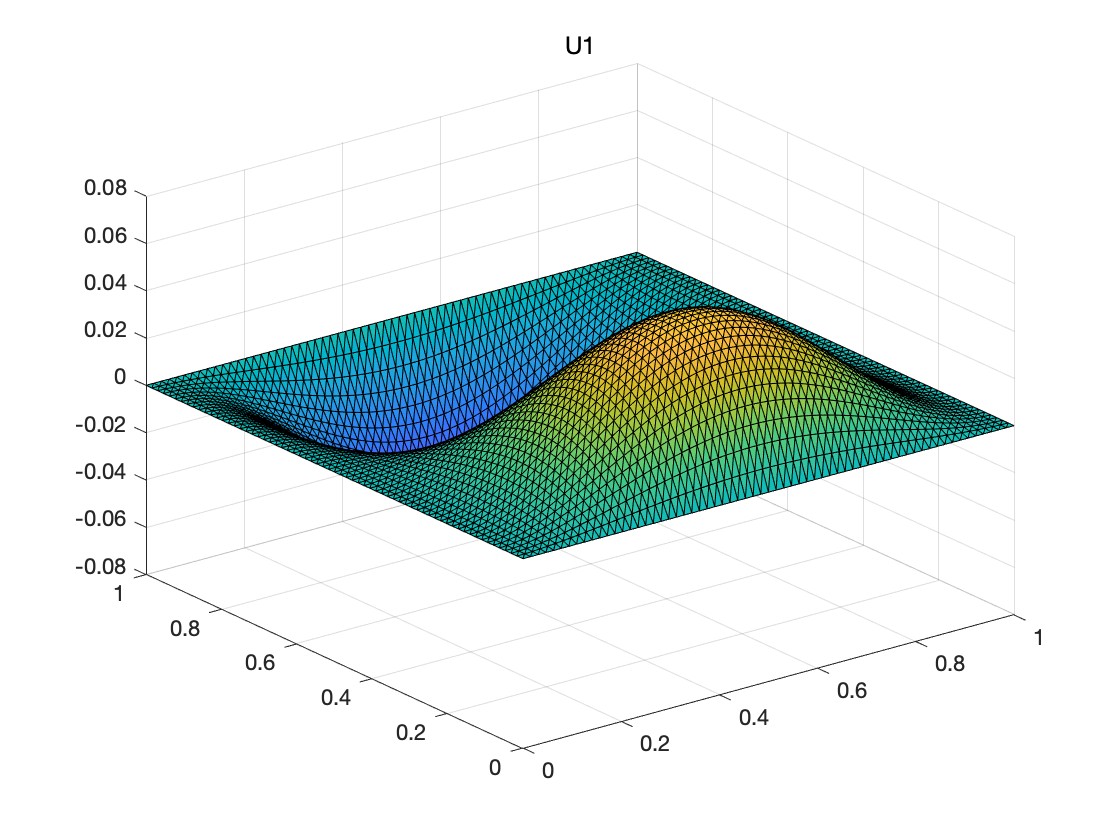}
        \caption{PR-EG approximation $u_{1,h}^{\mathrm{PR}}$.}
        \label{fig:preg_u1_num}
    \end{subfigure}

    \vskip\baselineskip

    \begin{subfigure}{0.45\textwidth}
        \centering
        \includegraphics[width=\textwidth]{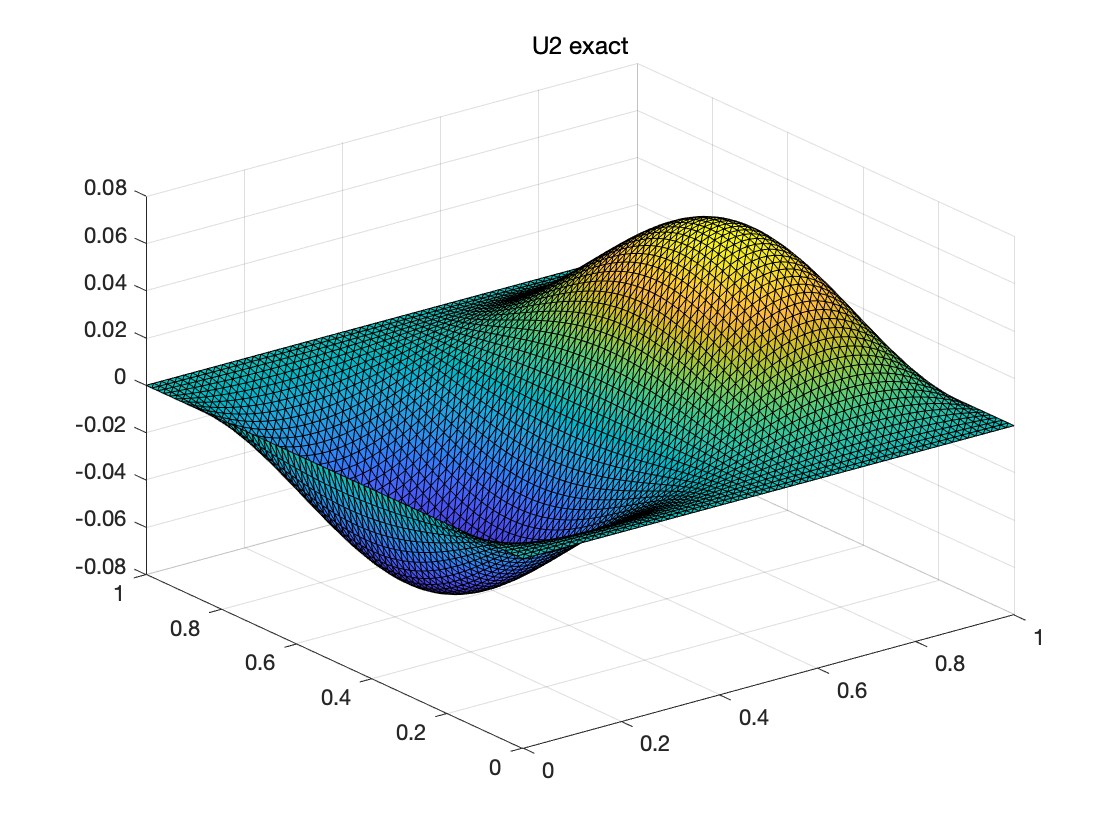}
        \caption{Exact velocity component $u_2$.}
        \label{fig:preg_u2_exact}
    \end{subfigure}
    \hfill
    \begin{subfigure}{0.45\textwidth}
        \centering
        \includegraphics[width=\textwidth]{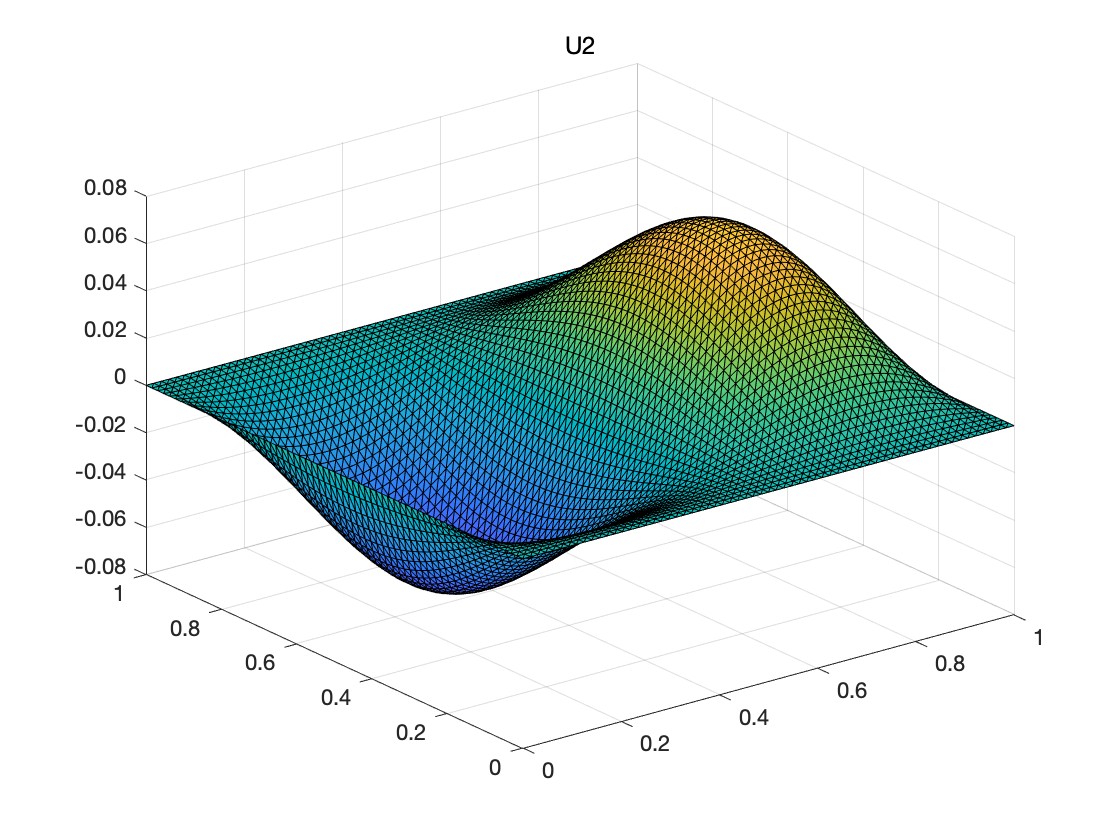}
        \caption{PR-EG approximation $u_{2,h}^{\mathrm{PR}}$.}
        \label{fig:preg_u2_num}
    \end{subfigure}

    \caption{Exact and PR-EG velocity components for Example~3.}
    \label{fig:preg_smooth_solution}
\end{figure}

For the PR-EG method, the numerical results are expected to exhibit
first-order convergence in the mesh-dependent energy norm,
\[
    \|\mathbf u-\mathbf u_h^{\mathrm{PR}}\|_{\mathcal E}=O(h),
\]
while the velocity $L^2$-error is expected to display the corresponding
second-order convergence behavior. The pressure approximation is also
expected to converge with first-order accuracy in the $L^2$-norm.\\
\indent More importantly, the comparison between the two methods is intended to
illustrate the effect of the velocity reconstruction operator. Whereas the velocity
accuracy of the standard EG method may be influenced by the pressure
approximation, the PR-EG velocity estimate contains no pressure-dependent
term. Therefore, the PR-EG curves are expected to exhibit a more uniform
behavior as $\mu$ decreases.\\
\indent Finally, Figure~\ref{fig:preg_smooth_solution} compares the exact velocity
components with the corresponding PR-EG approximations on a representative
mesh.

\subsubsection{Example 4: No-Flow Test with a Gradient Force}

The second example is designed to assess the pressure robustness of the
proposed PR-EG method. We consider a no-flow problem on the unit square
\[
    \Omega=(0,1)^2,
\]
for which the exact velocity is identically zero,
\[
    \mathbf{u}(x,y)=\mathbf{0},
\]
while the pressure is prescribed by
\[
    p_{\alpha}(x,y)
    =
    \frac{\alpha}{3}
    \left(
        x^3+y^3-\frac{1}{2}
    \right),
\]
where $\alpha>0$ is a scaling parameter controlling the magnitude of the
pressure. Since
\[
    \int_{\Omega}
    \left(
        x^3+y^3-\frac{1}{2}
    \right)\,dx\,dy=0,
\]
we have $p_{\alpha}\in L_0^2(\Omega)$.\\
\indent Because the exact velocity vanishes, both the viscous and nonlinear convective
terms are identically zero. Therefore, the corresponding forcing term is
given entirely by the pressure gradient,
\[
    \mathbf{f}_{\alpha}
    =
    \nabla p_{\alpha}
    =
    \alpha
    \begin{pmatrix}
        x^2\\
        y^2
    \end{pmatrix}.
\]
Hence, regardless of the magnitude of $\alpha$, the exact velocity remains
\[
    \mathbf{u}=\mathbf{0}.
\]

This example provides a direct test of pressure robustness. Indeed, the
gradient force $\mathbf{f}_{\alpha}=\nabla p_{\alpha}$ should be absorbed by
the pressure and should not generate a nonzero velocity field. However, for a
non-pressure-robust discretization, the discrete velocity may be contaminated
by the approximation of the pressure gradient, and this effect can become
particularly pronounced in the small-viscosity regime. In contrast, the
velocity reconstruction employed in the PR-EG formulation is designed to
eliminate the influence of the gradient component of the forcing term on the
discrete velocity.\\
\indent To examine this behavior, we fix
\[
    \mu=10^{-6},
    \qquad
    h=\frac{1}{64},
\]
and consider the pressure scaling parameters
\[
    \alpha
    =
    1,\;10^2,\;10^4,\;10^6.
\]
For each value of $\alpha$, the standard EG and PR-EG methods are solved on
the same mesh using the same numerical parameters and nonlinear stopping
criterion.\\
\indent The velocity errors obtained by the standard EG and PR-EG methods are
reported in Table~\ref{tab:no_flow}. A pressure-robust discretization is
expected to maintain a very small discrete velocity even when the magnitude
of the pressure gradient is increased by several orders of magnitude.

\begin{table}[htbp]
    \centering
    \caption{Comparison of the velocity differences of the standard EG and
    PR-EG methods for the no-flow test with $\mu=10^{-6}$ and $h=1/64$.}
    \label{tab:no_flow}
    \begin{tabular}{c|cc|cc}
        \hline
        $\alpha$
        & $\|\mathbf{u}_{h,\alpha}^{\mathrm{EG}}
        -\mathbf{u}\|_{\mathcal E}$
        & $\|\mathbf{u}_{h,\alpha}^{\mathrm{PR}}
        -\mathbf{u}\|_{\mathcal E}$
        & $\|\mathbf{u}_{h,\alpha}^{\mathrm{EG}}
        -\mathbf{u}\|_0$
        & $\|\mathbf{u}_{h,\alpha}^{\mathrm{PR}}
        -\mathbf{u}\|_0$
        \\
        \hline
        $10^{0}$
        & $4.93184\mathrm{e}{+00}$
        & $1.43669\mathrm{e}{-09}$
        & $1.56234\mathrm{e}{-01}$
        & $1.48411\mathrm{e}{-11}$ \\

        $10^{2}$
        & $1.65014\mathrm{e}{+01}$
        & $1.43671\mathrm{e}{-09}$
        & $4.55642\mathrm{e}{-01}$
        & $1.48419\mathrm{e}{-11}$ \\

        $10^{4}$
        & $8.45409\mathrm{e}{+02}$
        & $1.62928\mathrm{e}{-09}$
        & $7.13215\mathrm{e}{+00}$
        & $1.53243\mathrm{e}{-11}$ \\

        $10^{6}$
        & $3.24909\mathrm{e}{+04}$
        & $7.67842\mathrm{e}{-08}$
        & $5.34854\mathrm{e}{+01}$
        & $1.00325\mathrm{e}{-10}$ \\
        \hline
    \end{tabular}
\end{table}


\subsubsection{Example 5: Lid-driven Cavity Flow for the PR-EG Method}

In this example, we consider the classical lid-driven cavity flow in the unit square
$\Omega=(0,1)^2$. The body force is set to zero. A unit tangential velocity
\[
    \mathbf{u}=(1,0)^T
\]
is prescribed on the top boundary, while homogeneous Dirichlet boundary
conditions are imposed on the remaining three sides. The pressure is normalized
by the zero-mean condition
\[
    \int_{\Omega}p_h\,dx=0.
\]

\begin{figure}[H]
    \centering
    \includegraphics[width=0.80\textwidth]{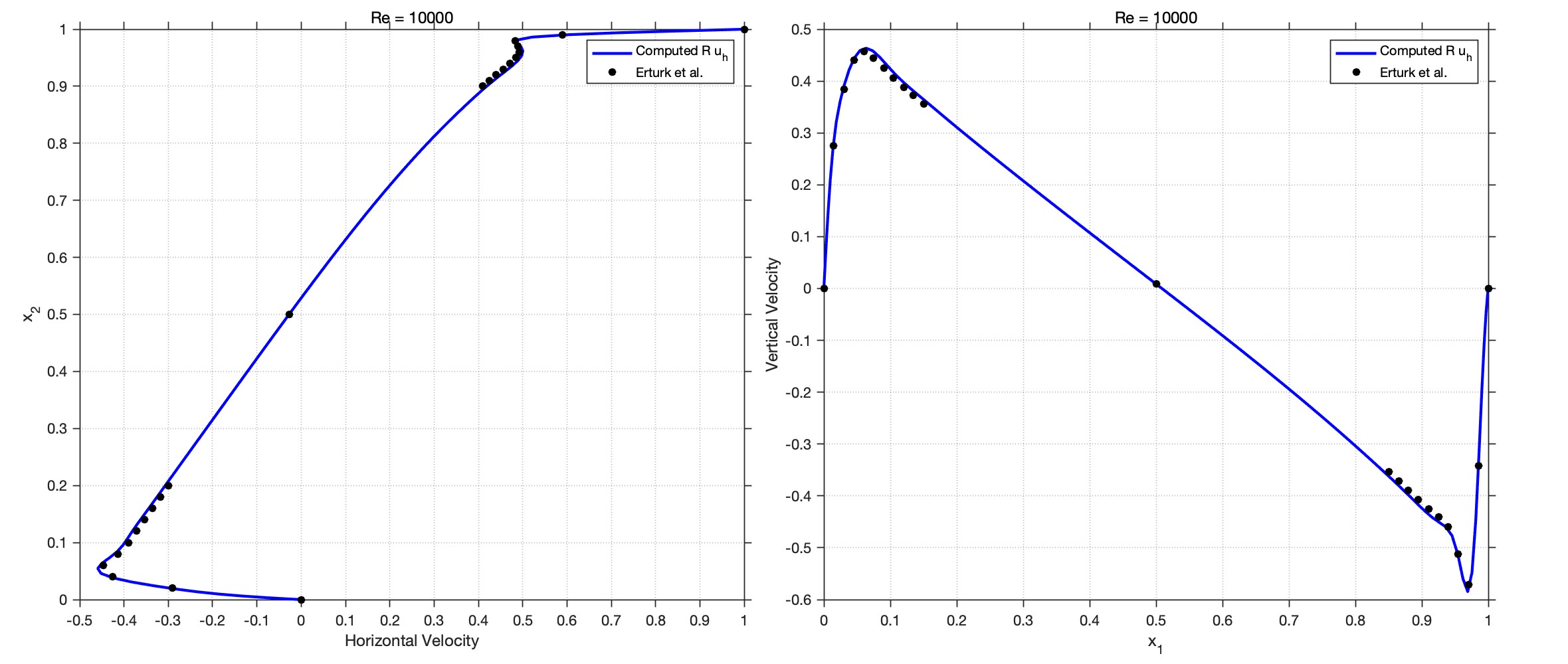}
    \caption{Centerline velocity profiles obtained by the PR-EG method for
    $\mu=10^{-4}$ ($\mathrm{Re}=10^{4}$). The reference data are taken from
    Erturk et al.~\cite{erturk2005numerical}.}
    \label{fig:preg_cavity_centerline_mu1e4}
\end{figure}

The PR-EG method is employed for the spatial discretization. The reconstruction
operator $\mathcal{R}$ is incorporated into the nonlinear convection term.
Since $\mathbf{f}=\mathbf{0}$, the reconstruction of the right-hand side does
not contribute explicitly in this example. A Stokes solution is first used as
the initial guess, and the resulting nonlinear system is solved by Newton
iteration.\\
\indent We first consider $\mu=10^{-4}$, corresponding to
$\mathrm{Re}=10^{4}$. The centerline velocity profiles are compared with the
reference data reported by Erturk et al.~\cite{erturk2005numerical}.
Figure~\ref{fig:preg_cavity_centerline_mu1e4} shows the comparison, demonstrating
good agreement between the PR-EG solution and the reference data.

\begin{figure}[H]
    \centering
    \includegraphics[width=0.55\textwidth]{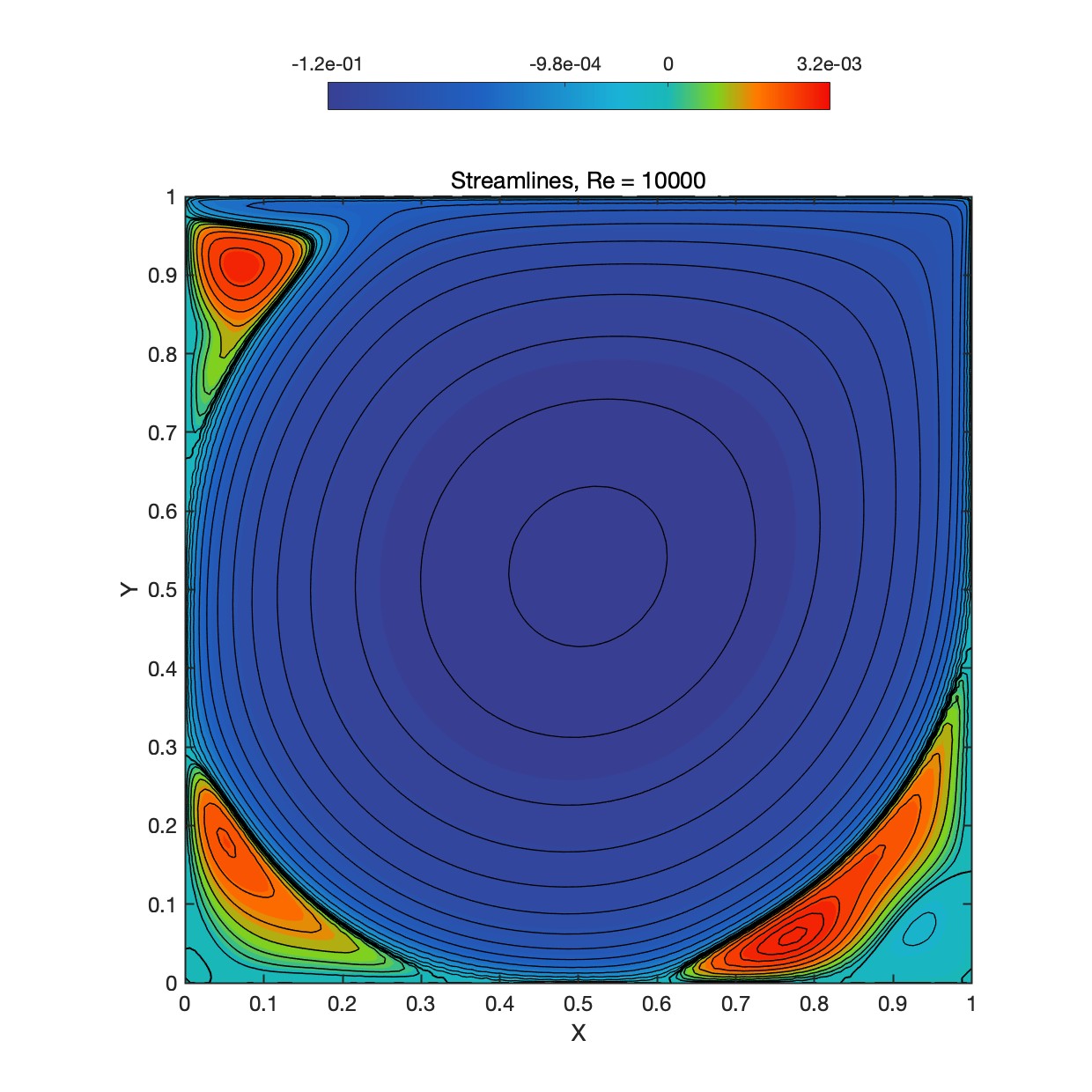}
    \caption{Streamline plot of the lid-driven cavity flow obtained by the
    PR-EG method for $\mu=10^{-4}$ ($\mathrm{Re}=10^{4}$).}
    \label{fig:preg_cavity_streamline_mu1e4}
\end{figure}

Figure~\ref{fig:preg_cavity_streamline_mu1e4} presents the corresponding
streamline plot for $\mu=10^{-4}$. The primary recirculation vortex and the
secondary corner vortices are clearly captured.

\begin{figure}[H]
    \centering
    \includegraphics[width=0.55\textwidth]{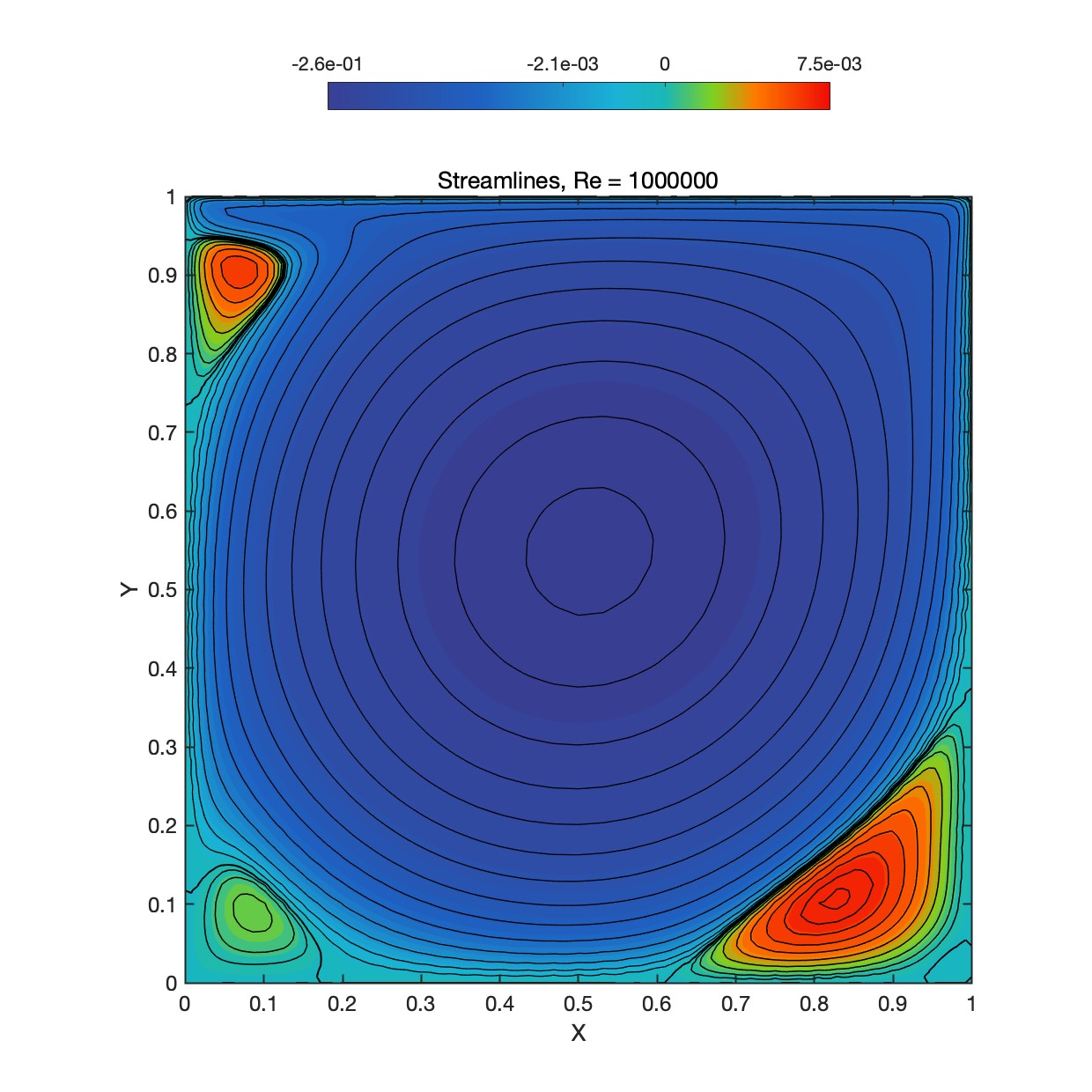}
    \caption{Streamline plot of the lid-driven cavity flow obtained by the
    PR-EG method for $\mu=10^{-6}$ ($\mathrm{Re}=10^{6}$).}
    \label{fig:preg_cavity_streamline_mu1e6}
\end{figure}

We further consider $\mu=10^{-6}$, corresponding to
$\mathrm{Re}=10^{6}$. Figure~\ref{fig:preg_cavity_streamline_mu1e6}
shows the resulting streamline plot. The numerical solution retains the
dominant recirculation structure and exhibits more pronounced flow structures
near the cavity corners.

\section{Conclusions}

In this work, we developed an enriched Galerkin method and a pressure-robust variant for the stationary incompressible Navier--Stokes equations.
The pressure-robust formulation employs a velocity reconstruction operator in the load term and nonlinear convection term, while the viscous and velocity--pressure coupling forms remain unchanged.
Under suitable regularity and small-data assumptions, stability and a priori error estimates are established.
In particular, the velocity error estimate of the PR-EG method is independent of the pressure approximation.
Numerical experiments confirm the expected convergence behavior and demonstrate the pressure robustness of the proposed method.
The lid-driven cavity tests further illustrate its applicability to convection-dominated flows.
Future work may consider time-dependent problems, higher-order discretizations, and more efficient nonlinear solvers.

    \section*{Data availability}
	Data will be made available on reasonable request.
	
	\section*{Declarations}
	The authors declare no competing interests.
	
	\bibliographystyle{unsrt}
	\bibliography{reference}
\end{document}